\documentclass[a4paper,11pt]{article}%,dvipdfm

\usepackage{mathrsfs} 

\usepackage{amsthm}
\usepackage{hyperref}
\usepackage{amssymb}
\usepackage{latexsym}
\usepackage{amsmath}
\usepackage{amsfonts}
\usepackage{mathabx}
\usepackage[dvips]{graphicx}
\usepackage[utf8]{inputenc}

\usepackage[T1]{fontenc}
\usepackage[francais,english]{babel}
\usepackage{url}
\usepackage{enumerate}
\usepackage{hyperref}
\usepackage{color}
\usepackage{bbm}

\usepackage{fancyhdr}
\newtheorem{Theorem}{Theorem}[part]
\newtheorem{Definition}{Definition}[part]
\newtheorem{Proposition}{Proposition}[part]

\newtheorem{Lemma}{Lemma}[part]
\newtheorem{Corollary}{Corollary}[part]
\newtheorem{Remark}{Remark}[part]

\makeatletter \@addtoreset{equation}{section}

\@addtoreset{Definition}{section}

\@addtoreset{Theorem}{section}

\@addtoreset{Proposition}{section}

\@addtoreset{Property}{section}

\@addtoreset{Assumption}{section}

\@addtoreset{Corollary}{section}

\@addtoreset{Lemma}{section}

\@addtoreset{Remark}{section}

\@addtoreset{Example}{section}

\newcommand{\cB}{\mathcal{B}}
\newcommand{\cC}{\mathcal{C}}
\newcommand{\cD}{\mathcal{D}}
\newcommand{\cE}{\mathcal{E}}
\newcommand{\cF}{\mathcal{F}}

\newcommand{\cI}{\mathcal{I}}

\newcommand{\cK}{\mathcal{K}}

\newcommand{\cM}{\mathcal{M}}

\newcommand{\cP}{\mathcal{P}}

\newcommand{\cR}{\mathcal{R}}

\newcommand{\cT}{\mathcal{T}}

\newcommand{\cV}{\mathcal{V}}

\newcommand{\cY}{\mathcal{Y}}
\newcommand{\cZ}{\mathcal{Z}}

\renewcommand{\P}{\mathbb{P}}

\newcommand{\R}{\mathbb{R}}

\newcommand{\dZ}
    {\ensuremath{\delta Z }}
    
\newcommand{\dY}
    {\ensuremath{\delta Y }}

\def \proof{{\noindent \bf Proof. }}
\def \eproof{\hbox{ }\hfill$\Box$}

\newcommand{\ud}{\mathrm{d}}

\newcommand{\1}{{\bf 1}}
    
\newcommand{\set}[1]
    {\ensuremath{\{ #1 \}}}
    
\newcommand{\HP}[1] %L2DPDT sur 0,T
    {\ensuremath{\mathscr{H}^{#1}}}    

\newcommand{\NH}[2] 
    { \ensuremath{ \left\|#2\right\|_{\mathscr{H}^{#1}} } }

\newcommand{\SP}[1]{\ensuremath{\mathscr{S}^{#1}}}

\newcommand{\NS}[2]
    {\ensuremath{\left\|#2\right\|_{ \mathscr{S}^{#1} }}}
\newcommand{\esp}[1]{\ensuremath{\mathbb{E} \!\! \left[#1\right] }}

\newcommand{\EFp}[2]
    {\ensuremath{% raccourci esperance conditionnelle discretisation
     \mathbb{E}_{#1}\!\!\left[#2\right] }}
     
\newcommand{\df}
    {\ensuremath{\delta f }}

\newcommand{\DY}
    {\ensuremath{\Delta Y }}
\newcommand{\DZ}
    {\ensuremath{\Delta Z }}

       \newcommand{\dxi}
    {\ensuremath{\delta \xi }}
       \newcommand{\dPsi}
    {\ensuremath{\delta \Psi }}

     \newcommand{\ixi}
    {\ensuremath{ {}^{i}\xi }}
    \newcommand{\exi}
    {\ensuremath{{}^{1}\!\xi }}
     \newcommand{\zxi}
    {\ensuremath{{}^{2}\!\xi }}
    \newcommand{\iY}
    {\ensuremath{ {}^{i}Y }}
    \newcommand{\eY}
    {\ensuremath{{}^{1}Y }}
     \newcommand{\zY}
    {\ensuremath{{}^{2}Y }}
    \newcommand{\iZ}
    {\ensuremath{ {}^{i}\!Z }}
    \newcommand{\eZ}
    {\ensuremath{{}^{1}\!Z }}
     \newcommand{\zZ}
    {\ensuremath{{}^{2}\!Z }}
        \newcommand{\iPhi}
    {\ensuremath{ {}^{i}\Phi}}
    \newcommand{\ePhi}
    {\ensuremath{{}^{1}\!\Phi }}
     \newcommand{\zPhi}
    {\ensuremath{{}^{2}\!\Phi }}    
           \newcommand{\iPsi}
    {\ensuremath{ {}^{i}\Psi}}
    \newcommand{\ePsi}
    {\ensuremath{{}^{1}\!\Psi }}
     \newcommand{\zPsi}
    {\ensuremath{{}^{2}\!\Psi }}   
     \newcommand{\ivf}
    {\ensuremath{ {}^{i}\!f }}
    \newcommand{\ef}
    {\ensuremath{ {}^{1}\!f }}
    \newcommand{\zf}
    {\ensuremath{ {}^{2}\!f }}
    
    \newcommand{\BP}[1]{\ensuremath{\mathscr{B}^{#1}}}
    
    \newcommand{\LP}[1]{\ensuremath{\mathscr{L}^{#1}}}
    
    \newcommand{\NB}[2]{\ensuremath{\left\| #2 \right\|_{\mathscr{B}^{#1}}} }

    \newcommand{\NL}[2]{\ensuremath{\left\| #2 \right\|_{\mathscr{L}^{#1}}} }

\newcommand{\pos}{\mathop{\mathrm{pos}}}
\DeclareMathOperator{\Span}{span}

\title{Obliquely Reflected Backward Stochastic Differential Equations}

\author{ 
          Jean-Fran\c{c}ois Chassagneux
         \\\small Laboratoire  de Probabilit\'{e}s, Statistique et Mod\'{e}lisation
         \\\small CNRS, UMR 8001,  Universit\'{e} Paris Diderot
         \\\small \sf jean-francois.chassagneux@univ-paris-diderot.fr
		\and
             Adrien Richou
             \\\small  \small Univ. Bordeaux, IMB, UMR 5251, F-33400 Talence, France.
              \\\small  \sf adrien.richou@math.univ-bordeaux.fr 
              }

\begin{document}
\maketitle

\begin{abstract}
In this paper, we study existence and uniqueness to multidimensional Reflected Backward Stochastic Differential Equations in an open convex domain, allowing for  oblique directions of reflection. In a Markovian framework, combining \emph{a priori} estimates for penalised equations and compactness arguments, we obtain existence results under quite weak assumptions on the driver of the BSDEs and the direction of reflection, which is allowed to depend on both $Y$ and $Z$. In a non Markovian framework, we obtain existence and uniqueness result for direction of reflection depending on time and $Y$. We make use in this case of stability estimates that require some smoothness conditions on the domain and the direction of reflection. 
\end{abstract}

\vspace{5mm}


\vspace{5mm}

\noindent {\bf MSC Classification (2000):}  93E20, 65C99, 60H30.

 %\newpage
 
 %\tableofcontents

%!TEX root = main.tex
\section{Introduction}

%\textcolor{magenta}{TODO pour la version à soumettre:
%\begin{enumerate}
%\item traiter les $\Phi_t \in \partial \varphi(Y_t)$...
%\end{enumerate}
%}

%Backward Stochastic Differential Equation have known an ever growing interest since the seminal paper by %\cite{Pardoux-Peng-90}. They were first introduced by Bismut in \cite{bis73} to solve stochastic control problem using a variational approach. Since then, the application to BSDE to the fields of stochastic control has revealed extremely fruitful. They are also a key tool in the study of some class of nonlinear parabolic PDEs and they can be used to design efficient numerical probabilistic method to approximate the solution of this PDEs. 

In this paper, we study a class of BSDE whose solution is constrained to stay in an open convex domain, hereafter denoted $\cD$. The ``reflection'' at the boundary of the domain is made along an oblique direction. Such equations are known as Obliquely Reflected BSDEs and they allow to represent the solution to some stochastic control problems. 
Precisely, let $(\Omega,\mathcal{F},\mathbb{P})$ be a complete probability space and $(W_t)_{t \in [0,T]}$ a $k$-dimensional Brownian motion, defined on this space, whose natural filtration is denoted $(\mathcal{F}_t)_{t \in [0,T]}$. $\mathcal{P}$ is the $\sigma$-algebra generated by the progressively measurable processes on $[0,T] \times \Omega$. In the sequel, $T>0$ is a terminal time for the equation under consideration. In this paper, we are interested in the study of existence and uniqueness of a $\mathcal{P}$- measurable solution $(Y,Z,\Phi)$ to the following equation
\begin{equation}
 \label{eq:EDSRORintro}
 \left\{ \begin{aligned} &(i)\;Y_t = \xi+\int_t^T f(s,Y_s,Z_s)\ud s-\int_t^T H(s,Y_s,Z_s)\Phi_s \ud s-\int_t^T Z_s \ud W_s, \; 0 \leqslant t \leqslant T,\\
         & (ii)\; Y_. \in \bar{\cD} \text{ a.s.}, \quad { \Phi_\cdot \in\partial \varphi(Y_\cdot) \;\ud t \otimes \ud \P-a.e.},\quad \int_0^T |\Phi_s|\1_{\set{Y_s \notin \partial \cD}} \ud s = 0,
        \end{aligned}
 \right.
\end{equation}
%
%\begin{equation}
% \label{eq GBSDE intro}
% %\left\{ \begin{aligned} &
% Y_t = \xi+\int_t^T f(s,Y_s,Z_s)ds-\int_t^T H(s,Y_s,Z_s)\Phi_s ds-\int_t^T Z_s dW_s, \quad 0 \leqslant t \leqslant T.%\\
%%         & Y_t \in \bar{\cD}, \quad \Phi_t \in\partial \varphi(Y_t),\quad \int_0^T |\Phi_t|\1_{\set{Y_t \notin \partial \cD}} \ud t = 0,
%%        \end{aligned}
% %\right.
%\end{equation}
where $\partial \cD$  is the boundary of the open convex domain $\cD$, $\varphi$ the (convex) indicator function of $\cD$, $\partial \varphi$ the subdifferential of $\varphi$ and $(f,H): \Omega \times [0,T] \times \mathbb{R}^d \times \mathbb{R}^{d \times k} \rightarrow (\mathbb{R}^d,\mathbb{R}^{d\times d})$ is a $\mathcal{P} \otimes \mathcal{B}(\mathbb{R}^d \times \mathbb{R}^{d \times k})$-measurable random function. The terminal value $\xi$ is given as a parameter and  belongs  to $\LP{2}(\cF_T)$, where for $p>0$ and a $\sigma$-algebra $\cB$,  $\LP{p}(\cB)$ is the space of $\cB$-measurable random $R$ variable satisfying $\esp{|R|^p}<+\infty$.
%\textcolor{red}{dimension + coefficient function}
Of course, we shall require some extra conditions to get an existence and uniqueness result. Classically, we will look for solution with the following  integrability property: $(Y,Z,\Phi) \in \SP{2}\times\HP{2}\times\HP{2}$, where,  for $p \in [1,\infty]$, $\HP{p}$ is the set of progressively measurable process $V$ such that
%\begin{align*}
$
\esp{(\int_0^T|V_t|^2 \ud t)^\frac{p}2} < +\infty\;,
$
%\end{align*}
and 
$\SP{p}$ is the set of continuous and adapted processes $U$ satisfying
%\begin{align*}
$
\esp{\sup_{t \in [0,T]} |U_t|^p }.
$
%\end{align*}
The main constraints on the couple $(Y,\Phi)$ are given in \eqref{eq:EDSRORintro}$(ii)$. As already mentioned, the first one is that $Y$ takes its value in $\bar{\cD}$, where $\cD$ is a non-empty open convex subset of $\R^d$. The fact that $\Phi_t \in\partial \varphi(Y_t)$ imposes that $\Phi$ is directed along the outward normal of the convex domain, the important point being that in $(i)$ this direction is perturbed by the operator $H$ and we are thus dealing with an oblique direction of reflection. When \eqref{eq:EDSRORintro}(i) is viewed backward in time, the process $\Phi$ or, more precisely $\Psi := H(\cdot)\Phi$, is the process allowing $Y$ to stay in $\bar{\cD}$. The condition  $\int_0^T |\Phi_t|\1_{\set{Y_t \notin \partial \cD}} \ud t = 0$ is then interpreted as a minimality condition, in the sense that $\Psi$ will be active only when $Y$ touches the boundary of the domain. This is of course one of the main ingredient to get uniqueness result for this kind of equation.

%\begin{align} \label{eq intro const 2}
%\Phi_t \in\partial \varphi(Y_t)\quad \text{and}\quad\int_0^T |\Phi_t|\1_{\set{Y_t \notin \partial \cD}} \ud t = 0,
%\end{align}
%where $\partial \varphi$ is the subdifferential operator of the convex function $\varphi$.
%The constraint in 
%\\
%Such kind of equations are coined Obliquely Reflected BSDE, due to the presence of matrix operator $H$ in \eqref{eq GBSDE intro}. 

Let us now mention some known results about these equations. In the one dimensional case, they have been first studied in \cite{el1997reflected} for the -- so called-- simply reflected case and in \cite{cvitanic1996backward} for the doubly reflected case. The literature on this specific form of equation has then grown very importantly due to their range of application, in particular in mathematical economics or stochastic control. The multidimensional case is only well understood in the case of normal reflection i.e. when the matrix-valued random function $H$ is equal to the identity, see \cite{Gegout-Petit-Pardoux-96}. The case of oblique direction of reflection has been only partially treated.  Up until recently, only very specific cases have been considered for the couple $(H,\cD)$. In \cite{ramasubramanian2002reflected}, the author studies the case of the reflection in an orthant with some restriction on the direction of oblique reflection and the driver $f$. Another case that has received a lot of 
attention is the setting of 
RBSDEs associated to ``switching problems'', see e.g. \cite{Hu-Tang-10, hamadene2010switching, Chassagneux-Elie-Kharroubi-11} and the references therein: the multidimensional domain has a specific form and the direction is along the axis, see also Section \ref{subse markov continuous} for more details. In this case, structural conditions on $f$ are required to retrieve existence and uniqueness results also. This restriction are based on the technique of proof used to obtain the results and which is mainly based on a monotonic limit theorem à la Peng \cite{peng1999monotonic}, in a multidimensional setting. To the best of our knowledge, the first attempt to treat the question of BSDEs with oblique reflection in full generality can be found in \cite{gassous2015multivalued}. Unfortunately, their setting is still quite restrictive concerning $H$ and $f$.

To the best of our knowledge, there is no, up to now, satisfying global approach for the question of well-posedness of Obliquely Reflected BSDEs, especially when compared the case of forward SDEs, where existence and uniqueness results are obtained for oblique reflection and general domain, see e.g. the seminal paper \cite{Lions-Sznitman-84}. 

Our goal in this paper is thus to prove  existence and uniqueness for the  RBSDEs \eqref{eq:EDSRORintro}
for generic $H$ and convex domain $\cD$ without imposing any structural dependence condition on the driver $f$ of the equation.
In this direction, we are able to obtain very general existence result in a Markovian setting, assuming only a weak domination property of the forward process, 
see Section \ref{se markov}. We also discuss there the non-uniqueness issue.
In the general case of $\cP$-measurable random coefficients $f$, $H$ and terminal condition $\xi$, we need to impose some smoothness assumptions on the domain, on $H$, which depend then only on the time and $y$ variables and on the terminal condition $\xi$.
In this case, we obtain both existence and uniqueness for the solution of \eqref{eq:EDSRORintro}. Let us remark that our new results on Obliquely Reflected BSDEs allow us to treat some new optimal switching problems called ``randomised switching problems'' and introduced in the Section 5 of \cite{Chassagneux-Richou-17}. 
\\
The main tool to obtain the existence result is to consider a sequence of penalised equation: for $n \in \mathbb{N}$, $t \in [0,T]$,
\begin{align}
 \label{eq:edsr penalized gene}
 Y^n_t &= \xi + \int_t^T f(s,Y^n_s,Z^n_s) \ud s -\int_t^T Z^n_s \ud W_s - \int_t^T H(s,Y^n_s,Z^n_s)\nabla \varphi^M_n (Y^n_s)ds\;,
\end{align}
where, for $y \in \R^d$, and some $M>0$,
 \begin{align}\label{de varphi n}
  \varphi^M_n(y) := n\inf_{x \in \mathcal{D}} \theta_M(y-x)\quad\text{ with }\quad
  \theta_M(h) =
\left \{
\begin{array}{rcl}
 M|h| {-\frac{M^2}2}& \text{if} & |h|>M,
 \\
{\frac12}|h|^2 & \text{if}  & |h| \le M.
\end{array}
\right.  
 \end{align}
The key point is to obtain the convergence in a strong sense of $(Y^n)$ to some process $Y$ along with some \emph{a priori} estimates on $(Z^n,\Phi^n)$. This will then allow to obtain the existence of some limiting process $(Z,\Phi)$ as well. In the setting of oblique direction of reflection, the question of uniqueness has generally to be investigated separately.   
\\
The first possible argument to obtain the convergence of $(Y^n)$ is to prove some monotonicity on the sequence to apply Peng's monotonic limit theorem \cite{peng1999monotonic}. In a multidimensional setting, this monotonicity is obtained under very restrictive structural condition on the coefficient. Nevertheless, it has been successfully used for the study of RBSDE associated to switching problem. Another possible argument is to invoke some fine compactness arguments and this is the approach followed in \cite{gassous2015multivalued}. But, again some strong structural conditions are required to obtain convergence results in a weak setting. In this paper, we follow a similar approach in the Markovian setting, see Section \ref{se markov}. At the heart of our proof, we use the paper \cite{Hamadene-Lepeltier-Peng-97}, which was concerned with multidimensional (non-reflected) BSDEs with continuous only driver $f$. With this approach, in the Markovian setting, we are able to obtain  existence result for $H$ that 
can depend on $Z$ and even be discontinuous. To the best of our knowledge, this is the first time such general setting is considered successfully. It has been brought to our attention that independantly from us, \cite{DeAngelis-Ferrari-Hamadene-17} has followed a similar approach to treat BSDEs associated to the classical switching problem in a more restrictive setting.

The last approach to obtain convergence of the sequence $(Y_n)$ is to show classically that it is a Cauchy sequence. This approach has been used in the case of  muldimensional RBSDE when there is no perturbation $H$ of the direction of reflection, namely $H$ is the identity matrix of $\R^d$, in the seminal paper \cite{Gegout-Petit-Pardoux-96}. To obtain this result and a key stability estimate, authors of \cite{Gegout-Petit-Pardoux-96} use dramatically the convexity property of the domain linked with the normal reflection by applying It\^o's formula to the Euclidean norm of the difference of two solutions. In our setting of general perturbation $H$, we cannot follow directly their proof. In order to retrieve the stability estimates, we modify the Euclidean norm to take into account the oblique reflection inspired by \cite{Lions-Sznitman-84}. Unfortunately, this produces new terms that have to be controlled. The most difficult one is certainly the term linked to the quadratic variation of the martingale term 
in \eqref{eq:EDSRORintro}(i) or \eqref{eq:edsr penalized gene}. Let us emphasize that this term cannot be dealt with as one would do in the forward SDE case. Nevertheless, we are able to treat this term using BMO martingales estimates. This tool was already used with success to deal with quadratic BSDEs but, to the best of our knowledge, this approach is completely new in the setting of Reflected BSDEs. We are then able to obtain in the non-Markovian setting existence results when $(\cD,H)$ satisfies some $C^2$ smoothness condition, with $H$ depending only on the time and $y$ variables. 
Let us note also that in this case the uniqueness result is obtained as an easy consequence of the stability estimate.

\vspace{5pt}
The rest of the paper is organised as follows. In the next Section, we present precisely our framework and the assumption made on the coefficients along with some discussions on these assumptions. We also prove the key \emph{a priori} and stability estimates, that will be used later on. In Section \ref{se smooth setting}, we present our first novel result on existence and uniqueness of Obliquely Reflected BSDEs in a regular setting for $(\cD,H)$. In Section \ref{se markov}, restricting to a Markovian framework, we extend our previous existence result assuming no regularity on $(\cD,H)$ and allowing a dependence in $Z$ for the $H$ operator. 
%Finally, in Section \ref{se the game you dont play}, we illustrate our result by applying them to a new class of randomised switching problem.

%\textcolor{magenta}{
%\begin{Remark}
%Equivalent characterisation if possible in terms of interior sphere?
%$C^2$ bounded domain satisfies the above condition...
%\end{Remark}
%}

%\vspace{5pt}
%In this case, we know that there exists a function $\phi \in \cC^2_b(\R^d,\R)$ that has the following properties
%\begin{align}
%&\cD = \set{x \in \R^d \,|\, \phi(x)<0},\; \partial \cD = \set{x \in \R^d \,|\, \phi(x)=0}
%\;,\; \R^d\setminus \bar{\cD} = \set{x \in \R^d \,|\, \phi(x) > 0} 
%\\
%&|\phi(x)| = d(x,\partial \cD) \quad \text{for} \quad x \in \cV \cup \R^d\setminus \bar{\cD} \;,\;\text{where $\cV$ is a neighbourhood of $\partial \cD$}
%%\\
%%
%\end{align}

\paragraph{Notations}
We denote by $\varphi$ the  indicator function of $\cD$
\begin{equation*}
 \varphi(y)= \begin{cases} 0 & \textrm{if } y \in \bar{\cD},\\
           +\infty & \textrm{otherwise,}
          \end{cases}
\end{equation*}
and $\partial \varphi$ its subdifferential operator:
\begin{equation*}
 \partial \varphi (y) = \begin{cases}
                         \left\{ \hat{y} \in \mathbb{R}^d : \hat{y} \cdot (z-y) \le 0, \forall z \in \bar{\cD} \right\} & \textrm{if } y \in \bar{\cD}\\
                        % \{0\}& \textrm{if } y \in \cD,\\
                         \varnothing & \textrm{if } y \notin \bar{\cD}.
                        \end{cases}
\end{equation*}
In particular, $\partial \varphi(y)$ is the closed  cone of outward normal to $\cD$ at $y$ when $y \in \partial \cD$ and $\partial \varphi(y)=\{0\}$ when $y \in \cD$.
Finally, we denote by $\mathfrak{P}$ the projection onto $\bar{\cD}$ and by $\mathfrak{n}(y)$ the set of unit outward normal at $y \in \partial \cD$.

For a matrix $M$, we denote $M^\dagger$ its transpose.

We denote $\BP{2}$, the set of processes $V \in \HP{2}$, such that% there exists a constant $C>0$ satisfying
\begin{align*}
 \NB{2}{V} := \NL{\infty}{\mathrm{sup}_{t \in [0,T]} \esp{\int_{t}^T |V_s|^2 \ud s | \cF_{t}}}^\frac12<+\infty.  \;
\end{align*}
Let us remark that $V \in \BP{2}$ means that the martingale $\int_0^. V_s \ud W_s$ is a BMO martingale and $\NB{2}{V}$ is the BMO norm of $\int_0^. V_s \ud W_s$. We refer to \cite{Kazamaki-94} for further details about BMO martingales.

%
%The set $\HP{p}$, $p \in [1,\infty]$ is the set of progressively measurable process $V$ such that
%\begin{align*}
%\NL{p}{(\int_0^T|V_t|^2 \ud t)^\frac12} < +\infty\;.
%\end{align*}

The set of continuous function from $[0,T]$ to $\R^n$ is denoted $\cC([0,T],\R^n)$. For $x \in \cC([0,T],\R^n)$, we denote by $\|x\|_\infty := \sup_{t \in [0,T]} |x_t|$, the sup-norm on this space.

%For a discrete grid $\pi := \set{t_0 := 0 < t_1 < \dots <t_m := T}$ of the interval $[0,T]$, we denote
%\begin{align}\label{eq notation disc path func}
%x^\pi := (x_{t_1}, \dots, x_{t_m}) \quad \text{ and } \quad x^\pi_t := (x_{t_1\wedge t}, \dots, x_{t_m}\wedge t)\;, 
%\end{align}
%for $0 \le t \le T$ and $x \in \cC([0,T],\R^n)$.
%\\
%We also introduce $\sharp \pi := m$ and $|\pi| := \max_{i<m}(t_{i+1}-t_i)$.
%
%
%\vspace{5mm}
%In the sequel $C$ is positive constant that may change from line to line and which depends implicitely on the parameter $b$, $T$ and $L$. We shall denote it $C_p$ if it depends on an extra parameter $p$.

%is essentially bounded.

%!TEX root = main.tex
\section{Setting and preliminary estimates}

In this section, we first introduce and discuss the main assumptions that will be used to obtain our existence and uniqueness results. In a second part, we give important \emph{a priori} estimates and prove a key stability result, which is one of the novelty in our approach to solve Obliquely Reflected BSDEs.

\subsection{Framework}
\label{subse framework}
The first minimal set of assumption that we consider here is the following.
\paragraph{Assumption (A)}
$ $
\begin{enumerate}[i)]
 \item $\xi$ is an $\mathcal{F}_T$-measurable random variable, $\mathbb{R}^d$-valued such that $\mathbb{E}\left[|\xi|^2\right]<+\infty$.
 \item $f: \Omega \times [0,T] \times \mathbb{R}^d \times \mathbb{R}^{d \times k} \rightarrow \mathbb{R}^d$ is a $\mathcal{P} \otimes \mathcal{B}(\mathbb{R}^d \times \mathbb{R}^{d \times k})$-measurable function and there exists a non negative progressively measurable process $\alpha \in \mathscr{H}^2(\mathbb{R})$ and a constant $L$ such that
 \begin{align} \label{eq bound f intro}
 |f(t,y,z)| \leqslant \alpha_t + L(|y|+|z|), \quad \forall (t,y,z) \in [0,T] \times \mathbb{R}^d \times \mathbb{R}^{d \times k}.
 \end{align}
 \item $H: \Omega \times [0,T] \times \R^d \times \mathbb{R}^{d \times k} \rightarrow \mathbb{R}^{d\times d}$ is a $\mathcal{P} \otimes \mathcal{B}(\mathbb{R}^d \times \mathbb{R}^{d \times k})$-measurable function and there exists a constant $\eta>0$ such that, for any $(t,y,z) \in [0,T] \times \mathbb{R}^d  \times \R^{d \times k}$
 %\begin{enumerate}[i)]
 %\item $H(t,y,z)$ is a symetric linear application,
 %\item 
 \begin{align}
  \label{inegalite H a intro}
 &  H(t,y,z)\upsilon \cdot \upsilon  \geq \eta , \quad \upsilon \in \mathfrak{n}(\mathfrak{P}(y))\;,%\quad \forall u \in \R^q, 
  \\
  &|H(t,y,z)| \leq L \label{inegalite bound H a intro}.
 \end{align}
 %\item $|H(t,\mathfrak{P}(y),z)| \leq L,$
% \item $H(t,y,z)=0$, when $y\in \cD$. \textcolor{red}{utile ?}
 %\end{enumerate}
%  \textcolor{red}{(on a pas besoin de symétrie a priori...)}, 
% \begin{equation}
%  \label{inegalite H a intro}
%  \langle H(t,\mathfrak{P}(y),z)(y-\mathfrak{P}(y)),y-\mathfrak{P}(y) \rangle \geq b |y-\mathfrak{P}(y)|^2,%\quad \forall u \in \R^q, 
% \end{equation}
 %\textcolor{red}{on en a pas besoin dans toutes les directions et on a pas besoin de l'inversibilite globale...}
% $$|H(t,\mathfrak{P}(y),z)| \leq L,$$
%\textcolor{magenta}{ $$H(t,y,z)=0, \quad \text{when } y\in \cD. ???$$}
\end{enumerate}

The above assumptions are too weak to obtain existence and uniqueness result in a general random framework. They will be used in Section \ref{se markov} in a Markovian framework  with their Markovian counterpart \textbf{(AM)}. Nevertheless, it is possible to derive useful \emph{a priori} estimates in the general setting of \textbf{(A)}.

\begin{Remark} \label{re extension H nonverysmooth}
%\begin{enumerate}[i)]
%\item \textcolor{magenta}{justification de la définition directe à tout l'espace malgré la donnée spécifique sur la frontière du domaine... => classification des convexes et tout ça... ???}
%%\item In view of the previous point, we define, for all $(t,y,z) \in [0,T] \times \mathbb{R}^d \times \mathbb{R}^{d\times k}$,
%%\begin{align}\label{eq de bar Hs}
%%H(t,y,z) = H(t,\mathfrak{P}(y),z).
%%\end{align}
%\end{enumerate}
In applications, $H(t,\cdot,z)$ is usually specified only on the boundary $\partial \cD$. The extension to $\R^d \setminus \bar{\cD}$ in a continuous way can be done easily by setting $H(t,y,z) := H(t,\mathfrak{P}(y),z)$. Moreover, if $H(t,\cdot,z)$ is a continuous and bounded function on $\partial \cD$ it is possible to extend it to a continuous and bounded function on $\bar{\cD}$. Indeed, $\bar{\cD}$ is homeomorph to a set $S$ which is a half plane of $\mathbb{R}^d$ or $\mathbb{R}^r \times B^{d -r}$, with $0 \leqslant r \leqslant d$. Moreover, the boundary of $\cD$ is sent to the boundary of $S$.  Then we  remark that the extension of $H(t,.,z)$ is straightforward when $\cD=S$. 
\end{Remark}

In the non-Markovian setting, our results require more smoothness and control on the parameters of the BSDE. We will then work under the following assumption.

\paragraph{Assumption (SB)} 
\begin{enumerate}[i)]
 \item $\xi$ is an $\mathcal{F}_T$-measurable $\bar{\cD}$-valued random variable and %belonging to 
 %$\LP{\infty}$.  \textcolor{red}{introduce class of admissible terminal condition here!}
the martingale $\cY^\xi_t := \EFp{t}{\xi}=\xi - \int_t^T\cZ^\xi_s \ud W_s$, $t \le T$, is BMO (see \cite{Kazamaki-94} for further details on BMO martingales). 
 
 \item $f: \Omega \times [0,T] \times \mathbb{R}^d \times \mathbb{R}^{d \times k} \rightarrow \mathbb{R}^d$ is a $\mathcal{P} \otimes \mathcal{B}(\mathbb{R}^d \times \mathbb{R}^{d \times k})$-measurable function, there exists  a constant $L>0$ such that,  for all $(t,y,y',z,z') \in [0,T] \times \mathbb{R}^d \times \mathbb{R}^d \times \mathbb{R}^{d \times k}\times \mathbb{R}^{d \times k}$,
 \begin{align} \label{eq hyp f snb}
 |f(t,y,z) - f(t,y',z')| &\le L \left(|y-y'| + |z-z'|\right)\;.
% \nonumber |f(t,y,z)| &\leqslant L(1+|z|).
 %\quad\text{and}
%\quad |f(t,0,0)| \leqslant \alpha_t, \quad %\forall (t,y,z) \in [0,T] \times \mathbb{R}^d \times \mathbb{R}^{d \times k}.
 \end{align}
Moreover, the process $\theta^\xi_\cdot := f(\cdot,\cY^\xi_\cdot,\cZ^\xi_\cdot)$ belongs to $\BP{2}$.

\item {The open convex domain $\cD$ is given by a $\cC^2(\R^d,\R)$ function $\phi$ with a bounded first derivative, namely  $\cD = \set{\phi<0}$ and $\partial \cD = \set{\phi = 0}$. $\phi$ is assumed . This function satisfies moreover %and $\partial \cD$ is $\cC^2$.
\begin{align}
|\phi(x)| = d(x,\partial \cD) \quad \text{for} \quad x \in \cV \cap (\R^d\setminus \bar{\cD}) \;,\;\text{where $\cV$ is a neighbourhood of $\partial \cD$.}
\end{align}
}
 \item $H:  [0,T] \times \R^d \rightarrow \mathbb{R}^{d\times d}$ is valued in the set of symmetric matrices $Q$ satisfying 
 \begin{align}\label{eq bound matrix H}
 |Q| \le L\;,\quad L |\upsilon|^2 \ge  \upsilon^\dagger Q \upsilon \ge \frac1L |\upsilon|^2\,,\;\forall \upsilon \in \R^d.
 \end{align}
 % and for some $\eta>0$, 
 %\begin{equation}
 % \label{inegalite H a intro smooth}
 % H(t,\mathfrak{P}(y))\upsilon \cdot \upsilon  \geq \eta,\quad \forall y \in \R^d, t \in [0,T], \upsilon \in \mathfrak{n}(\mathfrak{P}(y)).
 %\end{equation}
 $(t,y) \mapsto H(t,y)$ is a $\cC^{0,1}$-function  and $(t,y) \mapsto H^{-1}(t,y)$ is a  $\cC^{1,2}$ function satisfying
 \begin{align}\label{eq ass bound deriv H}
 |\partial_y H| + |H^{-1}| + |\partial_t H^{-1}| + |\partial_y H^{-1}| + |\partial^2_{yy} H^{-1}| \le \Lambda\,,
 \end{align}
 for some positive $\Lambda$.
 %\begin{enumerate}[i)]
 %\item $H(t,y,z)$ is a symetric linear application,
 %\item 
 %\item 

% \item $H(\cdot)$ is a symmetric and invertible and
%  \begin{align}\label{eq ass bound H}
%\quad  \text{ and }  \quad
% |\partial_t H^{-1}| + |\partial_y H^{-1}| + |\partial^2_{yy} H^{-1}| \le L\,.
% \end{align}
% 
 %\end{enumerate}
%\item Assumption $(A)$ holds and moreover, the $\bar{\cD}$-valued  random variable $\xi$ in $(A.1)$ belongs to $\LP{\infty}(\cF_T)$, the processus $\alpha$ in $(A.2)$ belongs to $\BP{2}$.

\end{enumerate}

We first comment assumptions made on BSDE parameters.
\begin{Remark}\label{re comment hyp xi}
% Let us observe that under the BMO condition, there exists $\mu^\xi > 0$, such that 
%$\esp{e^{\mu^\xi \sup_{t \in [0,T]}|\cY^\xi_t|} } < \infty$ and that $\NB{2}{\cZ^\xi} <\infty$. Moreover, the process $\theta^\xi_\cdot := f(\cdot,\cY^\xi_\cdot,\cZ^\xi_\cdot)$ belongs to $\BP{2}$.
%For later use, we define
 %\begin{align}\label{eq de sup Yxi expo int}
%\sigma^\xi := \esp{e^{\mu^\xi \sup_{t \in [0,T]}|\cY^\xi_t|} } + \NB{2}{\cZ^\xi} + \NB{2}{\theta^\xi}< \infty \;.
 %\end{align}
%\textcolor{blue}{ %a virer puis a remettre:
\begin{enumerate}[i)]
\item  Let us observe that under the BMO condition, there exists $\mu^\xi > 0$, such that 
$\esp{e^{\mu^\xi \sup_{t \in [0,T]}|\cY^\xi_t|} } < \infty$ and that $\NB{2}{\cZ^\xi} <\infty$, see e.g. \cite{Kazamaki-94}.
For later use, we define
 \begin{align}\label{eq de sup Yxi expo int}
\sigma^\xi := \esp{e^{\mu^\xi \sup_{t \in [0,T]}|\cY^\xi_t|} } + \NB{2}{\cZ^\xi} + \NB{2}{\theta^\xi}< \infty \;.
 \end{align}
\item Condition \textbf{(SB)}(ii) is a mix of the property of $\xi$, $f$ and the domain $\cD$. In many applications, it will be straightforward to check. For example, it is trivially satisfied in the following cases:
\begin{enumerate}
\item $\sup_{y \in \cD}f(s,y,z) \le C$, for some $C>0$;
\item $\xi \in \LP{\infty}(\cF_T)$;
\item $\cD$ is a bounded domain\;.
\end{enumerate}
%\item \textcolor{red}{moment exp $(\cY^\xi)^*$}
\item If \textbf{(SB)} holds, then \textbf{(A)} holds as well. Indeed, one can set
%\begin{align*}
$
\alpha := L\left( |\cY^\xi| + |\cZ^\xi| + |\theta^\xi| \right) \;.
$
%\end{align*}
\end{enumerate}
%}
\end{Remark}

We now discuss the various assumptions made on $H$ and the domain $\cD$.
\begin{Remark}
\label{rem hyp domain}
%\textcolor{blue}{Plus besoin depreciser la definition fonction $C^2$ sur la variete}
\begin{enumerate}[i)]
\item The function $\phi$ can be constructed as in e.g. \cite{Gegout-Petit-Pardoux-96}  Section 2.4 if the convex domain $\cD$ is $\cC^2$. From $(SB)(iii)$, it follows that $\partial \phi (x)$ (resp. $\mathfrak{n}(x)$) is the outward normal (resp. unit outward normal) of $\cD$  at a point $x \in \partial \cD$. Moreover, since $\cD$ is convex, $\phi$ is convex on $\R^d\setminus\cD$ and thus, $\partial_{xx}^2\phi$ is a positive semi-definite matrix  on this domain. Let us also observe that the application $\mathfrak{P}: \R^d \setminus \bar{\cD}\rightarrow \partial \cD$ is $\cC^2$.
%\item \textcolor{red}{$\bar{H}$ to define is smooth ? to define at all?}
\item The matrix $H$ defines on $\partial \cD$ a  unit vector field $\nu$ in the following way
\begin{align*}
\tilde{\nu}(t,y) := H(t,y) \mathfrak{n}(y)  \;\text{ and }\; \nu(t,y) := \frac{\tilde{\nu}(t,y)}{|\tilde{\nu}(t,y)|}\;,\quad \text{ for } y \in \partial \cD\;,
\end{align*}
which represents the oblique direction of reflection. Then, \eqref{inegalite H a intro} rewrites as
\begin{align}
\langle \tilde{\nu}(t,y), \mathfrak{n}(y) \rangle \ge \eta\;,\quad \text{ for } y \in \partial \cD\;. \label{eq vector field property}
\end{align}
 In applications, it is generally the case that only the smooth vector field $\nu$ is given on $\partial \cD$. Following \cite{Lions-Sznitman-84}, it is possible to construct $H$ satisfying $(SB)(iv)$ on $\partial \cD$ and then to extend it on $\bar{\cD}$ under $(SB)(iii)$ using classical extension results, see e.g. \cite{gildbarg1977elliptic}.
 %Lemma of \textcolor{red}{Gildbarg-Trudinger}.
\end{enumerate}
\end{Remark}

We now introduce a class of terminal conditions that are admissible for the purpose of our work, in the sense that we can obtain an existence and uniqueness result for this class.

\begin{Definition}\label{de new class}
For $\beta >0$, the class $\mathfrak{T}_\beta$ is the subset of $\xi \in \LP{2}(\cF_T)$ satisfying: there exists $\lambda_\xi >  \beta$, such that
\begin{align}\label{eq mini exp integrability on Z}
 \esp{e^{\lambda_\xi \int_0^T |\cZ^{\xi}_s|^2 \ud s}} < \infty \;,
\end{align}
where $\cZ^\xi$ is given by the martingale representation theorem applied to $\cY^\xi_t := \EFp{t}{\xi}=\xi - \int_t^T\cZ^\xi_s \ud W_s$, $t \le T$.
\end{Definition}

We study the class $\mathfrak{T}_\beta$ in Section \ref{subse term xi}. Especially, we exhibit some specific elements of this class that are quite useful for applications.

%\textcolor{red}{FAIRE une remarque sur les constantes, quid de la dependance en L???}

\vspace{10pt}
\begin{Remark}\label{re about the constant} 
{
In the following, we will use in proofs the notation ``$C$'' to denote a generic constant that may change from line to line and that depends in an implicit way on $T$, $L$ and $\eta$. We shall denote it $C_\theta$, if it depends on an extra parameters $\theta$. In the statement of the results, we prefer the notation ``$c$'' and the dependence upon any extra parameters on top of $T$, $L$ and  $\eta$, will also be made clear.
}
\end{Remark}
%%%%%
%%%%%

\subsection{\emph{A priori} estimates}
\label{subse a priori estimates}
In this section, we prove some \emph{a priori} control on the solution  $(Y,Z,\Phi)\in \SP{2}\times\HP{2}\times \HP{2}$ to the following generalised BSDE
 \begin{align}
 \label{eq a priori bounds}
 Y_t = \xi+\int_t^T f(s,Y_s,Z_s)ds-\int_t^T H(s,Y_s,Z_s)\Phi_s ds-\int_t^T Z_s dW_s, \quad 0 \leqslant t \leqslant T\;.
 %,\\
        % & Y_t \in \bar{\cD}, \quad \Phi_t \in\partial \varphi(Y_t),\quad |\Phi_t|\1_{\set{Y_t \notin \partial \cD}} = 0,
  %      .
\end{align}
Importantly, we assume that $(Y,Z,\Phi)$ satisfies the following structural condition:
 \begin{align} \label{eq hyp gen Phi}
 \EFp{t}{\int_t^T|\Phi_s|^2 \ud s} \le K \EFp{t}{\int_t^T|f(s,Y_s,Z_s)|^2\ud s} \;, \text{ for some } K>0\;.
 \end{align}
 
 \vspace{4pt}
Equation \eqref{eq a priori bounds} encompasses both the obliquely reflected BSDE \eqref{eq:EDSRORintro} and its penalised approximation given in equation \eqref{eq:edsr penalized gene}. The key point for these two equations will then be to prove that their solutions satisfy condition \eqref{eq hyp gen Phi}.

\vspace{5pt} 
Our first estimate is quite classical. 
\begin{Lemma} \label{le first basic a priori}
Assume \textbf{(A)}. Let $(Y,Z,\Phi)\in \SP{2}\times\HP{2}\times \HP{2}$ be a solution to \eqref{eq a priori bounds} with condition \eqref{eq hyp gen Phi}  holding true.
 Then, for some $c:=c(K)$,
 \begin{align*}
 \sup_{t\in[0,T]}  \esp{ |Y_t|^2} + \esp{\int_0^T|Z_s|^2\ud s } \le c \,\esp{|\xi|^2+\int_0^T |\alpha_s|^2 ds}\;.
 \end{align*}
\end{Lemma}

\proof
We apply It\^o's formula to $|Y|^2$ to obtain
\begin{align}
|Y_t|^2 +  \int_t^T \!\!\! |Z_s|^2 \ud s 
&= |\xi|^2 +
2\int_t^T  \!\!\! Y_s f(s,Y_s,Z_s) \ud s 
%\\
%&
-2 \int_t^T  \!\!\! Y_s H(s,Y_s,Z_s)\Phi_s \ud s - 2\int_t^T  \!\!\! Y_s Z_s \ud W_s.
\end{align}
We observe, thanks to the square integrability of $Y$ and $Z$ that $\int_0^\cdot Y_s Z_s \ud W_s$ is a true martingale. This yields, %for all $r \le t$,
\begin{align*}%\label{eq temp level 0}
\EFp{}{|Y_t|^2 +  \int_t^T \!\!\! |Z_s|^2 \ud s }
&=
\EFp{}{|\xi|^2 + 2\int_t^T  \!\!\! Y_s f(s,Y_s,Z_s) \ud s 
- 2 \int_t^T  \!\!\! Y_s H(s,Y_s,Z_s)\Phi_s \ud s }\;.
\end{align*}
\\
We thus compute, using \eqref{eq hyp gen Phi} and Assumption \textbf{(A)}(ii), the boundedness of $H$ and Young's inequality, for some $ \epsilon \in (0,1)$,
\begin{align*}
\EFp{}{\int_t^T  \!\!\! Y_s H(s,Y_s,Z_s)\Phi_s \ud s} 
&\le
C_{K}\EFp{}{ \int_t^T(\frac1\epsilon |Y_s|^2 + \epsilon|\Phi_s|^2 ) \ud s  }
\\
&\le
C_{K}\EFp{}{ \int_t^T( \frac1\epsilon |Y_s|^2 + \epsilon|Z_s|^2  + |\alpha_s|^2) \ud s  }\;.
\end{align*}
Similarly, we get
\begin{align*}
\EFp{}{ \int_t^T  \!\!\! Y_s f(s,Y_s,Z_s) \ud s } \le  C\EFp{}{ \int_t^T(\frac1\epsilon |Y_s|^2 + \epsilon|Z_s|^2  + |\alpha_s|^2) \ud s  }\;.
\end{align*}
For $\epsilon$ small enough and using Gronwall Lemma, we deduce
\begin{align*}
\EFp{}{|Y_t|^2 +  \int_t^T \!\!\! |Z_s|^2 \ud s } \le 
C_{K} \EFp{}{|\xi|^2 + \int_t^T|\alpha_s|^2 \ud s}\;.
\end{align*}
\eproof

\vspace{10pt}

The following proposition refines the previous estimates in the smooth setting of Assumption \textbf{(SB)}. It will also allow to use the stability result proved in the next section. Interestingly, it shows that  most of the properties of the martingale $\cY^\xi$ are transferred to the non-linear process given in equation \eqref{eq a priori bounds}.

\begin{Proposition}\label{pr second a priori estimate}
Assume that \textbf{(SB)} holds. Let $(Y,Z,\Phi)\in \SP{2}\times\HP{2}\times \HP{2}$ be a solution to \eqref{eq a priori bounds} with condition \eqref{eq hyp gen Phi}  in force. Then,  the following holds
 \begin{enumerate}[i)]
 \item $(Y,Z,\Phi) \in \SP{2}\times\BP{2}\times \BP{2}$ with, for some 
 $c:=c(K,\sigma^\xi)$,
  \begin{align}\label{eq prop 2.1 easy}
 \esp{e^{\mu^\xi \sup_{t \in [0,T]}|Y_t|} } +
 \NB{2}{\Phi} + \NB{2}{Z} \le c\;,
 \end{align}
 and,  for all $b>0$ and some $c':=c'(b,K,\sigma^\xi)$
  \begin{align}\label{eq easy peasy}
  \esp{e^{b \int_0^T \left( |\Phi_s| + |Z_t-\cZ_t^\xi| + |\theta^\xi_t| \right) \ud s}} \le c' \;.
 \end{align}
 \item  Moreover, if $\xi \in \mathfrak{T}_\beta$, for some $\beta>0$, then there exists $\Theta \in \HP{2}$ such that,  for all non negative increasing process $\gamma$ satisfying $\esp{|\gamma_T|^p} < \infty$ for some $p > 1$ (depending on $\gamma$), we have for all $t \in [0,T]$
 \begin{align}\label{eq pr bien vu}
  \EFp{t}{\int_t^T \gamma_s |Z_s|^2 \ud s} \le \EFp{t}{\int_t^T \gamma_s |\Theta_s| \ud s} < +\infty 
 \end{align}
 and  for some $\lambda \in ( \beta, \lambda^\xi)$ and $c:=c(K,\sigma^\xi,\lambda)$,
 \begin{align}\label{eq control exp}
\esp{e^{\lambda \int_0^T |\Theta_t| \ud t }} \le c \;.
 \end{align}
 \end{enumerate}
\end{Proposition}

%\begin{bluetext}
\proof An important step to obtain our estimates below is to compare the BSDE $(Y,Z)$ with the martingale $\cY^\xi$. To this end, we introduce for this proof $\DY := Y- \cY^\xi$ and $\DZ = Z - \cZ^\xi$. 
%For the reader's convenience, we also denote $\kappa := \NB{2}{\theta^\xi}$ and $\ell:=%\NB{2}{\cZ^\xi}$\\
\\1.a We apply It\^o's formula to $|\DY|^2$ to obtain
\begin{align*}
|\DY_t|^2 +  \int_t^T \!\!\! |\DZ_s|^2 \ud s 
&=
2\int_t^T  \!\!\! \DY_s f(s,Y_s,Z_s) \ud s 
%\\
%&
-2 \int_t^T  \!\!\! \DY_s H(s,Y_s)\Phi_s \ud s - 2\int_t^T  \!\!\! \DY_s \DZ_s \ud W_s.
\end{align*}
We observe, thanks to the square integrability of $\DY$ and $\DZ$ that $\int_0^\cdot\DY_s \DZ_s \ud W_s$ is a true martingale. This yields, for all $r \le t$,
\begin{align}\label{eq temp level 0}
\EFp{r}{|\DY_t|^2 +  \int_t^T \!\!\! |\DZ_s|^2 \ud s }
&=
2\EFp{r}{\int_t^T  \!\!\! \DY_s f(s,Y_s,Z_s) \ud s 
- \int_t^T  \!\!\! \DY_s H(s,Y_s)\Phi_s \ud s }
\end{align}
\\
We thus compute, using \eqref{eq hyp gen Phi}, the Lipschitz continuity of $f$, the boundedness of $H$ and Young's inequality, for all $r\le t$ and some $\epsilon \in (0,1)$
\begin{align*}
\EFp{r}{\int_t^T  \!\!\! \DY_s H(s,Y_s)\Phi_s \ud s} 
&\le
C_{K}\EFp{r}{ \int_t^T(\frac1\epsilon |\DY_s|^2 + \epsilon|\Phi_s|^2 ) \ud s  },
\\
&\le
C_{K}\EFp{r}{ \int_t^T(\frac1\epsilon |\DY_s|^2 + \epsilon|\DZ_s|^2  + |\theta^\xi_s|^2) \ud s  }\;.
\end{align*}
Similarly, we obtain
\begin{align*}
\EFp{r}{\int_t^T  \!\!\! \DY_s f(s,Y_s,Z_s) \ud s }
\le
C\EFp{r}{ \int_t^T( \frac1\epsilon |\DY_s|^2 + \epsilon|\DZ_s|^2  + |\theta^\xi_s|^2) \ud s  } \;.
\end{align*}
Combining the last two estimates with \eqref{eq temp level 0}, setting $\epsilon$ small enough and using Gronwall Lemma, we get for all $r\leqslant t$
\begin{align*}%\label{eq control DY DZ}
\EFp{r}{|\DY_t|^2 + \frac12 \int_t^T \!\!\! |\DZ_s|^2 \ud s }
&\le
C_{K}\EFp{r}{\int_t^T|\theta^\xi_s|^2 \ud s }.
\end{align*}
1.b  Setting $r=t$ in the previous inequality, we have
\begin{align}\label{eq control DY DZ}
\sup_{t \in [0,T]} |\DY_t|^2 + \NB{2}{\DZ}^2 \le C_{K,\sigma^\xi} %\NB{2}{\theta^\xi}^2 \;,
\end{align}
from which we straightforwardly deduce
\begin{align}
\esp{e^{\mu^\xi\sup_{t\in[0,T]}|Y_t|}} &\le 
\esp{e^{\mu^\xi\sup_{t\in[0,T]}(|\DY_t| + |\cY^\xi_t|)}} \le C_{K,\sigma^\xi}
\\
\text{ and }\;\NB{2}{Z} &\le \NB{2}{\cZ^\xi} +  \NB{2}{\DZ} \le C_{K,\sigma^\xi}\;.\label{eq control Z}
\end{align}
Combining \eqref{eq control DY DZ}  with \eqref{eq hyp gen Phi}, we obtain
\begin{align}\label{eq control Phi B2 norm}
\NB{2}{\Phi} \le C_{K,\sigma^\xi}\;.
\end{align}
This concludes the proof of \eqref{eq prop 2.1 easy}.\\
2.a
\textcolor{black}{We denote $\mathfrak{R}:= |\Phi| + |\DZ|+|\theta^\xi|$. 
For all $b>0$, we use Young inequality to get
\begin{align}
\label{eq pre john}
 \esp{e^{b \int_0^T\mathfrak{R}_s \ud s}}
 \le e^{\frac{b^2T}{\varepsilon}}\esp{e^{\varepsilon \int_0^T |\mathfrak{R}_s|^2 \ud s}}
\end{align}
for all $\varepsilon >0$. Then, by setting $\varepsilon =\left( 1+4\NB{2}{\Phi}^2+ 4\NB{2}{\DZ}^2 + 4\NB{2}{\theta^\xi}^2\right)^{-1}$ we compute, for all $r \in [0,T]$,
$$\EFp{r}{\int_r^T \varepsilon |\mathfrak{R}_s|^2 \ud s} \le 3 \varepsilon  \EFp{r}{ \int_r^T |\Phi_s|^2 + |\DZ_s|^2 + |\theta^\xi_s|^2 \ud s} \leqslant \frac{3}{4}.$$
Going back to \eqref{eq pre john}
 and applying the John-Nirenberg formula, see Theorem 2.2 in \cite{Kazamaki-94},
 we obtain
 \begin{align}\label{eq theta is done}
 \esp{e^{b \int_0^T \mathfrak{R}_s \ud s}}
 \le C_{K,\sigma^\xi,b}\;,
 \end{align}
 which proves \eqref{eq easy peasy}.
}\\
\noindent 2.b Applying It\^o's formula to $\gamma_\cdot|\DY_\cdot|^2$, on $[t,T]$, we compute,
\begin{align}\label{eq ref gamma}
\gamma_t|\DY_t|^2 +  \int_t^T\gamma_s |\DZ_s|^2 \ud s 
+ \int_t^T|\DY_s|^2 \ud \gamma_s
=&
2\int_t^T \gamma_s\DY_s f(s,Y_s,Z_s) \ud s 
\\
&-2 \int_t^T \gamma_s\DY_s H(s,Y_s)\Phi_s \ud s - 2\int_t^T \gamma_s\DY_s \DZ_s \ud W_s. \nonumber
\end{align}
Let us observe that the local martingale $\int_0^{\cdot}\gamma_t\DY_t \DZ_t \ud W_t$ is a true martingale. 
%\textcolor{magenta}{On a pas besoin de BDG}
Indeed, we compute, using Burkholder-Davis-Gundy inequality,
\begin{align*}
\esp{\sup_{s \in [0,T]} |\int_0^{s}\gamma_t\DY_t \DZ_t \ud W_t|}
&\le C \esp{(\int_0^{T}|\gamma_t\DY_t \DZ_t|^2 \ud t)^\frac12}
\\
&\le C_{K,\sigma^\xi} \esp{|\gamma_T|  \left(\int_0^{T}|\DZ_t|^2 \ud t\right)^\frac12}
\end{align*}
where we used \eqref{eq control DY DZ} for the last inequality.
Using  H\"older inequality, denoting $q$ the conjugate exponent of $p$, we get
\begin{align*}
\esp{\sup_{s \in [0,T]} |\int_0^{s}\gamma_t\DY_t \DZ_t \ud W_t|} &
\le C_{K,\sigma^\xi,p} \esp{|\gamma_T|^p}^\frac1p \esp{ \left(\int_0^{T}|\DZ_t|^2 \ud t\right)^{\frac{q}2}}^\frac1q \;.
\end{align*}
From the energy inequality, c.f. (VI.109.7) in \cite{Dellacherie-Meyer-85}, we have that
\begin{align*}
\esp{ \left(\int_0^{T}|\DZ_t|^2 \ud t\right)^{\frac{\lceil q \rceil}2}} \le C_{q} \NB{2}{\DZ}^{\lceil q \rceil} \;.
\end{align*}
We thus deduce
\begin{align*}
\esp{\sup_{s \in [0,T]} |\int_0^{s}\gamma_t\DY_t \DZ_t \ud W_t|}
< \infty \;.
\end{align*}
%thanks to the boundedness of $\dY$ and the square integrability property of $\gamma$ %and $\dZ$.
Since $\gamma$ is non-decreasing, we then compute, using \eqref{eq ref gamma}, \eqref{eq hyp gen Phi} and the Lipschitz continuity of $f$,
\begin{align}\label{eq majo by XiZ}
\EFp{t}{\int_t^T\gamma_s |\DZ_s|^2 \ud s}
\le
\EFp{t}{
\int_t^T \gamma_s \Xi_s \ud s 
} < \infty \;,
\end{align}
where we set
%\begin{align*}
\textcolor{black}{
$
\Xi := C_{K,\sigma^\xi}(1 + \mathfrak{R})
$
%\end{align*}
recalling that $\DY$ is bounded by \eqref{eq control DY DZ} and $\mathfrak{R}$ is defined in step 2.a.  
%\textcolor{blue}{
%It follows from step 2.a that, for all $b>0$,}
%\begin{align*}
%\esp{e^{b\int_0^T|\DZ_s| \ud s}} \le C_{L,K,\sigma^\xi,b} \; \quad \text{and} \quad \esp{e^{b\int_0^T|\theta^\xi_s| \ud s}} \le C_{L,K,\sigma^\xi,b}\;.
%\end{align*}
%Subsequently, combining the previous inequality with 
Using \eqref{eq theta is done}
%and using Cauchy-Schwarz inequality, 
we compute
\begin{align}\label{eq exp control any XiZ}
\esp{e^{b\int_0^T|\Xi_s| \ud s}} \le C_{K,\sigma^\xi,b}  \;,
\end{align}
for all $b > 0$.
}
\\
2.c We set $\lambda = (1+\epsilon)\beta$ with $\epsilon > 0$ such that $(1+\epsilon)^2\beta \le \lambda_\xi$, recalling Definition \ref{de new class}. Now we  define
 %small enough such that, 
 $$ \Theta := (1+\epsilon)|\cZ^\xi|^2 + (1+\frac1\epsilon)\Xi.$$ 
We observe that $  \EFp{t}{\int_t^T\gamma_s |\Theta_s|^2 \ud s} < \infty $: this follows from \eqref{eq majo by XiZ} and the fact that
 \begin{align*}
  \EFp{t}{\int_t^T\gamma_s |\cZ^\xi_s|^2 \ud s} < \infty \;.
   \end{align*}
This last inequality is simply obtained by applying It\^o's formula to $\gamma |\cY^\xi|^2$ which yields
\begin{align*}
\esp{\int_0^T\gamma_s |\cZ^\xi_s|^2 \ud s} \le \esp{\gamma_T|\xi|^2}
\le \NL{p}{\gamma_T} \NL{q}{|\xi|^2} <\infty\;.
\end{align*}
From the definition of $\Theta$, we have that, for $t \le T$,
 $$ \EFp{t}{\int_t^T \gamma_s \Theta_s \ud s} 
 \ge  (1+\epsilon)\EFp{t}{\int_t^T \gamma_s |\cZ^\xi_s|^2 \ud s} 
 + (1+\frac1\epsilon) \EFp{t}{\int_t^T \gamma_s |\DZ_s|^2 \ud s} $$
 where we used \eqref{eq majo by XiZ}. Then it follows from Young's inequality,
 $$
   \EFp{t}{\int_t^T \gamma_s |Z_s|^2 \ud s} \le 
   C\EFp{t}{\int_t^T \gamma_s \Theta_s \ud s} 
   < \infty \,,$$
   which proves \eqref{eq pr bien vu}.
 Finally, we compute using H\"older's inequality,
 \begin{align*}
  \esp{e^{(1+\epsilon)\beta \int_0^T\Theta_s \ud s}} \le C\esp{e^{ (1+\epsilon)^2\beta\int_0^T|\cZ^\xi_s|^2 \ud s} }^\frac1{1+\epsilon}\esp{e^{(1+\frac1{\epsilon})^2\beta\int_0^T|\Xi_s|\ud s}}^\frac{\epsilon}{1+\epsilon}
 \end{align*}
and using the fact that $\xi \in \mathfrak{T}_\beta$ and \eqref{eq exp control any XiZ}, we obtain \eqref{eq control exp} with $\lambda = (1+\epsilon)\beta$.
 %\begin{align}\label{eq control Theta Z}
  %\esp{ e^{(1+\epsilon)\beta \int_0^T |\Theta^Z_s| \ud s}} \le C_{K,\sigma^\xi,(1+\varepsilon)\beta}\;.
 %\end{align}
\eproof
%\end{bluetext}

\vspace{5pt}
Let us remark the following result, that will be useful in the next section.

\begin{Corollary}\label{co Y bounded}
Assume that \textbf{(SB)} holds. Let $(Y,Z,\Phi)\in \SP{2}\times\HP{2}\times \HP{2}$ be a solution to \eqref{eq a priori bounds} with condition \eqref{eq hyp gen Phi}  in force and assume moreover that $\xi \in \LP{\infty}(\cF_T)$ then $Y$ is bounded, namely for some $c:=c(K,\sigma^\xi,\NL{\infty}{\xi})$, we have
\begin{align*}
\sup_{t \in [0,T]} |Y_t| \le c\;.
\end{align*}
\end{Corollary}
\proof We observe that $|\cY^\xi_t| \le \NL{\infty}{\xi}$ and then conclude using \eqref{eq control DY DZ}.
\eproof

\subsection{A stability result}

In this section, we prove a key estimate for the difference of two solutions of the
generalised BSDE \eqref{eq a priori bounds} satisfying \eqref{eq hyp gen Phi}. For $i \in \set{1,2}$, we denote $(\iY,\iZ,\iPhi)$ the solutions associated to parameters $(\ixi,\ivf)$ and we furthermore assume that
\begin{align}\label{eq prop of iPhi}
\iPhi_. \in \partial \varphi(\mathfrak{P}(\iY_.))\; \ud \P\otimes\ud t -\text{a.e.} \quad\text{ and }\quad \int_0^T|\iPhi_t|\1_{\set{\iY_t \in \cD}} \ud t = 0\;.
\end{align}

\begin{Remark} The above assumption allows us to  cover both cases of equation \eqref{eq:EDSRORintro} and equation \eqref{eq:edsr penalized gene}.
\end{Remark}

\vspace{4pt}

We now define $\dY = \eY-\zY$, $\dZ = \eZ-\zZ$, $\dPsi = \ePsi - \zPsi$, where $\iPsi = H(\cdot,\iY)\iPhi$ and $\df = \ef(\cdot,\eY,\eZ) - \zf(\cdot,\zY,\zZ)$. We have the following key result for our work.

%\textcolor{red}{Def $\Lambda$}.

\begin{Proposition} \label{pr stability} %\textcolor{red}{ecrire pour $(\theta,\tau)$ a la place de $(g,T)$}
Assume that $\textbf{(SB)}$ holds. There exist two increasing functions $\mathfrak{B}(\cdot)$
and $\mathfrak{A}(\cdot)$ from $(0,\infty)$ to $(0,\infty)$,
such that  
%for some $\alpha \ge \mathfrak{A}(\Lambda)$
%and $\beta \ge \mathfrak{B}(\Lambda)$, $t \in [0,T]$, 
for all $\exi$ belonging to $\mathfrak{T}_{\mathfrak{B}(\Lambda)}$,
setting
{
\begin{align*}
\Gamma_t := e^{ \mathfrak{A}(\Lambda) t + \int_{0}^t \mathfrak{B}(\Lambda)\set{ \Theta^{1}_s + \Theta^\Phi_s +\Theta^f_s}\ud s}
\end{align*}
where $\Theta^\Phi := |\ePhi| + |\zPhi|$,  $\Theta^f := |\eZ-\cZ^{\exi}| + |\theta^{\exi}|$, and $\Theta^1$ is given by Proposition \ref{pr second a priori estimate}(ii) applied to the BSDE with parameters $(\exi,\ef)$,
}
%\begin{align*}
%\Gamma_t := e^{ \mathfrak{A}(\Lambda) t + \int_{0}^t \mathfrak{B}(\Lambda)%\set{ \Theta^{\eZ}_s + \Theta^\Phi_s }\ud s}
%\end{align*}
%
%with $\Theta^\Phi := |\ePhi| + |\zPhi|$
we have,
\begin{enumerate}[i)]
\item  $\esp{|\Gamma_T|^p}< c$ for some $p := p(\Lambda) > 1$ and $c:= c(K,\Lambda,\sigma^{\exi},\sigma^{\zxi});$
\item for some $c':=c'(K,\Lambda,\sigma^{\exi},\sigma^{\zxi})$, and for all $t \le T$, 
%\textcolor{blue}{
\begin{align}\label{eq pr stab}
 |\dY_t|^2 + \EFp{t}{\int_t^T |\dZ_s|^2 \ud s} \le& c'
\EFp{t}{\Gamma_T|\dxi|^2
+
\int_t^T \Gamma_s |(\ef-\zf)(s,\eY_s,\eZ_s)|^2 \ud s
\right. 
\\
& \left. + \int_t^T \!\!\Gamma_s \left( |\mathfrak{P}(\zY_s) - \zY_s|+|\mathfrak{P}(\eY_s) - \eY_s|)
(|\ePhi_s| +  |\zPhi_s| \right)
\ud s
}. \nonumber
\end{align}
\end{enumerate}

%}
%\textcolor{red}{estime conditionnelle => T2.2 Kazamaki pour appli}
\end{Proposition}

%\begin{bluetext}
\proof
%\textcolor{red}{minoration ci dessous vrai après avoir pris esp conditionnelle}
% we set, for $ t \le T$,
% \begin{align*}
% \Gamma_t = e^{ \alpha t + \int_{T-\delta}^t \beta\set{ |\ePhi_s| + |\zPhi_s|  + |\eZ_s|^2 }\ud s}
% \end{align*}
% for $\alpha$ to be fixed later on \textcolor{red}{depending only on $\Lambda$}. 
% Let us observe then, that from \eqref{hy:momentexpo},
% \begin{align}\label{eq bound Gamma2}
% \esp{|\Gamma_T|^2} \le C_{\Lambda,L}\;.
% \end{align}
In this proof, we denote $A = (a^{ij}) = H^{-1}$ and the following simplified notation will be used $a^{ij}_t = a^{ij}(t,\eY_t)$, $\partial_t a^{ij}_t = \partial_ta^{ij}(t,\eY_t)$, $\partial_y a^{ij}_t = \partial_ya^{ij}(t,\eY_t)$,
$\partial^2_{yy} a^{ij}_t = \partial^2_{yy}a^{ij}(t,\eY_t)$ and $f_t = f(t,\eY_t,\eZ_t)$. For the reader's convenience, we shall also denote $\sigma := \sigma^{\exi} \vee \sigma^{\zxi}$ in the proof below.
\\
1. We first show the integrability property of $\Gamma$. We first recall that from
Proposition \ref{pr second a priori estimate}, for all $b>0$, we have
\begin{align}\label{eq let's get started}
\esp{e^{b\int_0^T\set{|\Theta^\Phi_s| + |\Theta^f_s|}\ud s}} \le C_{K,\sigma,b} \;.
\end{align}
Setting $p:=p(\Lambda)> 1$ such that $p^2\mathfrak{B}(\Lambda)\le \lambda_\xi$, recall Definition \ref{de new class}, we obtain using H\"older inequality,
\begin{align*}
\esp{|\Gamma_T|^p} 
&\le C_\Lambda \esp{e^{p\mathfrak{B}(\Lambda)\int_0^T \set{ \Theta^{1}_s + \Theta^\Phi_s + \Theta^f_s}\ud s }} 
\\
&\le C_{K,\Lambda,\sigma} \esp{e^{p^2\mathfrak{B}(\Lambda)\int_0^T \Theta^{1}_s \ud s }}^\frac1p 
\\
&\le
C_{K,\Lambda,\sigma} \,,
\end{align*}

%
%\begin{align}\label{eq let's get started}
%\esp{e^{b\int_0^T|\Theta^\Phi_s| \ud s}} \le C_{L,K,\sigma,b} \;.
%\end{align}
%Setting $p:=p(\Lambda)> 1$ such that $p^2\mathfrak{B}(\Lambda)\le \lambda_\xi$, recall Definition \ref{de new class}, we obtain using H\"older inequality,
%\begin{align}
%\esp{|\Gamma_T|^p} 
%&\le C_\Lambda \esp{e^{p\mathfrak{B}(\Lambda)\int_0^T \set{ \Theta^{\eZ}_s + \Theta^\Phi_s }\ud s }} 
%\\
%&\le C_{L,K,\Lambda,\sigma} \esp{e^{p^2\mathfrak{B}(\Lambda)\int_0^T \Theta^{\eZ}_s \ud s }}^\frac1p 
%\\
%&\le
%C_{L,K,\Lambda,\sigma} \,,
%\end{align}
where we used \eqref{eq let's get started} and Proposition \ref{pr second a priori estimate} (ii)\;.

\noindent 2.a To obtain the stability result, we first expand the product $(\Gamma_t \delta Y_t^\dagger A(t,Y_t) \delta Y_t)_{0 \le t \le T}$.
Applying It\^o's formula, we compute, for $1 \le i,j \le d$,
\begin{align}
&\ud[\Gamma_t a^{ij}_t\dY^i_t\dY^j_t] /\Gamma_t
\\ \nonumber
&=
\dY^i_t\dY^j_t\left( a^{ij}_t \mathfrak{A}(\Lambda) + \partial_ta^{ij}_t \right) \ud t =:(\cE^\cT_t)^{ij} \ud t 
\\ \nonumber
&+ \dY^i_t\dY^j_t\left( a^{ij}_t \mathfrak{B}(\Lambda) \Theta^1_t+ \frac12 \mathrm{Tr}[\partial^2_{yy}a^{ij}_t \eZ_t\eZ^\dagger_t ]\right) \ud t
=: (\cE^Z_t)^{ij} \ud t 
\\ 
&+\set{a^{ij}_t(-\dY^i_t\df^j_t-\dY^j_t\df^i_t + \mathfrak{B}(\Lambda)\Theta^f_t\dY^i_t\dY^j_t)  
-\partial_y a^{ij}_t \ef_t \dY^i_t\dY^j_t}\ud t =:(\cE^f_t)^{ij} \ud t  \label{eq f term}
\\
&+ \set{a^{ij}_t\sum_{m=1}^k\dZ^{im}_t\dZ^{jm}_t + \sum_{m=1}^k\partial_y a^{ij}_t \eZ^{.m}_t(\dY^i_t\dZ^{jm}_t+\dY^j_t\dZ^{im}_t)}\ud t =: (\cE^{\dZ}_t)^{ij} \ud t  \label{eq dZ term}
%\\
%&+\frac12\dY^i_t\dY^j_t\mathrm{Tr}[\partial^2_{yy}a^{ij}_t \eZ_t\eZ^\dagger_t ]  \ud t
\\
&+\set{a^{ij}_t(\dY^i_t\dZ^{j.}_t+\dY_t\dZ^{i.}_t) + 
\dY^i_t\dY^j_t\partial_y a^{ij}_t \eZ_t}\ud W_t  =: \ud M^{ij}_t \label{eq mart term}
\\
&+\set{ a^{ij}_t(\dY^i_t\dPsi^j_t+\dY^j_t\dPsi^i_t) 
+(\partial_y a^{ij}_t \ePsi_t  +  \mathfrak{B}(\Lambda) a^{ij}_t\Theta^\Phi_t )\dY^i_t\dY^j_t} \ud t =:(\cE^\cR_t)^{ij} \ud t \;.\label{eq reflect term}
\end{align}
%\begin{align}
%&\ud[\Gamma_t a^{ij}_t\dY^i_t\dY^j_t] /\Gamma_t
%\\
%&=
%\dY^i_t\dY^j_t\left( a^{ij}_t \mathfrak{A}(\Lambda) + \partial_ta^{ij}_t \right) \ud t =:(\cE^\cT_t)^{ij} \ud t 
%\\
%&+ a^{ij}_t\dY^i_t\dY^j_t\left(\mathfrak{B}(\Lambda) \Theta^Z_t+ \frac12 \mathrm{Tr}[\partial^2_{yy}a^{ij}_t \eZ_t\eZ^\dagger_t ]\right) \ud t
%=: (\cE^Z_t)^{ij} \ud t 
%\\
%&+\set{a^{ij}_t(-\dY^i_t\df^j_t-\dY^j_t\df^i_t )  
%-\partial_y a^{ij}_t\ef_t \dY^i_t\dY^j_t}\ud t =:(\cE^f_t)^{ij} \ud t  \label{eq f term}
%\\
%&+ \set{a^{ij}_t\sum_{k=1}^m\dZ^{ik}_t\dZ^{jk}_t + \sum_{k=1}^m\partial_y a^{ij}_t \eZ^{.k}_t(\dY^i_t\dZ^{jk}_t+\dY^j_t\dZ^{ik}_t)}\ud t =: (\cE^{\dZ}_t)^{ij} \ud t  \label{eq dZ term}
%%\\
%%&+\frac12\dY^i_t\dY^j_t\mathrm{Tr}[\partial^2_{yy}a^{ij}_t \eZ_t\eZ^\dagger_t ]  \ud t
%\\
%&+\set{a^{ij}_t(\dY^i_t\dZ^{j.}_t+\dY_t\dZ^{i.}_t) + 
%\dY^i_t\dY^j_t\partial_y a^{ij}_t \eZ_t}\ud W_t  =: \ud M^{ij}_t \label{eq mart term}
%\\
%&+\set{ a^{ij}_t(\dY^i_t\dPsi^j_t+\dY^j_t\dPsi^i_t) 
%+(\partial_y a^{ij}_t \ePsi_t  + \mathfrak{B}(\lambda) a^{ij}_t\Theta^\Phi_t )\dY^i_t\dY^j_t} \ud t =:(\cE^\cR_t)^{ij} \ud t \;.\label{eq reflect term}
%\end{align}
We now study each term separately.

\noindent 2.b We start by the reflection terms in \eqref{eq reflect term}.
We first observe that
\begin{align*}
\sum_{1\le i,j \le d} (\cE^\cR_t)^{ij} =  2 A(t,\eY_t) \dY_t \cdot \dPsi_t  &+ \!\! \sum_{1\le i,j \le d} \partial_y a^{ij}_t \ePsi_t 
+ \mathfrak{B}(\Lambda) \Theta^\Phi_t \dY_t \cdot A(t,\eY_t)  \dY_t\;.
\end{align*}
Recalling \eqref{eq bound matrix H} and \eqref{eq ass bound deriv H}, we compute
\begin{align}
\sum_{1\le i,j \le d} (\cE^\cR_t)^{ij} \ge 2 A(t,\eY_t) \dY_t \cdot \dPsi_t + \left(\frac{\mathfrak{B}(\Lambda)}{L}\Theta^\Phi_t - C_\Lambda|\ePhi_t| \right)  |\dY_t|^2\;. \label{eq reflection term base}
\end{align}

For the first term in the right hand side of \eqref{eq reflection term base},  we compute
\begin{align}  
A(t,\eY_t) \dY_t \cdot \dPsi_t  &= A(t,\eY_t) \dY_t \cdot \ePsi_t - A(t,\zY_t) \dY_t \cdot \zPsi_t -\set{A(t,\eY_t)- A(t,\zY_t)}\dY_t \cdot \zPsi_t .\label{eq reflection term base 2}
\end{align}
We now observe that, %using the symmetry of $A$ and \eqref{eq prop convex}
\begin{align}%\label{eq reflection term part 1}
A(t,\eY_t) \dY_t \cdot \ePsi_t &= \dY_t \cdot A(t,\eY_t) \ePsi_t 
= \dY_t \cdot  \ePhi_t \nonumber\\
& \ge ( \eY_t - \mathfrak{P}(\eY_t)+ \mathfrak{P}(\zY_t) - \zY_t) \cdot  \ePhi_t \nonumber
\\
& \ge - \left(  |\mathfrak{P}(\eY_t) - \eY_t| + |\mathfrak{P}(\zY_t) - \zY_t|\right)|\ePhi_t| \label{eq reflection term first}
\end{align}
where we used \eqref{eq prop of iPhi} and the convexity property of $\cD$.
Similarly, we compute
\begin{align}
-A(t,\zY_t) \dY_t \cdot \zPsi_t  \ge - \left(  |\mathfrak{P}(\eY_t) - \eY_t| + |\mathfrak{P}(\zY_t) - \zY_t| \right)|\zPhi_t|. \label{eq reflection term second}
\end{align}
For the last term in the right-hand side of \eqref{eq reflection term base 2}, we get, using the Lipschitz property of $A$ that
\begin{align}
\set{A(t,\eY_t)- A(t,\zY_t)}\dY_t \cdot \zPsi_t \ge - C_\Lambda |\dY_t|^2 |\zPhi_t|. \label{eq reflection term third}
\end{align}
Combining \eqref{eq reflection term first}-\eqref{eq reflection term second}-\eqref{eq reflection term third}
with  \eqref{eq reflection term base}, we obtain, %using \eqref{eq useful upper bound}, 
for $\mathfrak{B}(\Lambda)$ large enough,
\begin{align}
\EFp{t}{\int_t^T \Gamma_s \!\!\!\sum_{1\le i,j \le d} (\cE^\cR_s)^{ij} \ud s}\ge  - C_\Lambda \EFp{t}{\int_t^T \!\!\Gamma_s(|\mathfrak{P}(\zY_s) - \zY_s|
+  |\mathfrak{P}(\eY_s) - \eY_s|) (|\ePhi_s| +|\zPhi_s|) \ud s }.\label{eq reflection term final}
\end{align}

2.c Using Young's inequality, we compute, recalling \eqref{eq bound matrix H} and \eqref{eq ass bound deriv H},
\begin{align} \label{eq lower bound EdZ}
\sum_{1 \le i,j \le d}(\cE^{\dZ}_t)^{ij} \ge 
\frac1{2L} |\dZ_t|^2 -C_\Lambda|\dY_t|^2|\eZ_t|^2 \;.
\end{align}
The terms $\cE^f$ in \eqref{eq f term} can be lower bounded, using Young's inequality and \eqref{eq bound matrix H}, by
{
\begin{align} \label{eq lower bound Ef}
\sum_{1 \le i,j \le d}(\cE^{f}_t)^{ij}  \ge & -\frac1{3L} |\dZ_t|^2 +\left(\frac{\mathfrak{B}(\Lambda)}L \Theta^f_t- C_{K,\Lambda,\sigma}(1+ |\theta^{\exi}_t|+|\eZ_t-\cZ^{\exi}|)\right)|\dY_t|^2 
\nonumber
\\
&- C_{\Lambda}|(\ef-\zf)(t,\eY_t,\eZ_t)|^2\;,
\end{align}
recalling \eqref{eq control DY DZ}.
}
%
%\begin{align} \label{eq lower bound Ef}
%\sum_{1 \le i,j \le d}(\cE^{f}_t)^{ij}  \ge -\frac1{3L} |\dZ_t|^2 -C_{\Lambda}(1+|\eZ_t|^2)|\dY_t|^2 - C|(\ef-\zf)(t,\eY_t,\eZ_t)|^2\;.
%\end{align}
%\textcolor{magenta}{ajouter un terme quand $f$ n'est pas bornee en $y$.} 
We also have that
\begin{align}\label{eq lower bound Ez and EcT}
 \sum_{1 \le i,j \le d}(\cE^{Z}_t)^{ij}  \ge \left( \frac{\mathfrak{B}(\Lambda)}{L}\Theta_t^1-C_\Lambda |\eZ_t|^2\right)|\dY_t|^2
%\end{align}
%and
%\begin{align}\label{eq lower bound EcT}
\;\text{ and }
\sum_{1 \le i,j \le d}(\cE^{\cT}_t)^{ij}  \ge \left( \frac{\mathfrak{A}(\Lambda)}{L}-C_\Lambda \right)|\dY_t|^2\;.
\end{align}
%\begin{bluetext}
2.d We now consider the local martingale $\cM$, defined by
\begin{align*}
 \cM_t := \sum_{1\le i,j \le d}  \cM^{ij}_t \; \text{ with } \; \cM^{ij}_t := \int_0^t \Gamma_s \ud M^{ij}_s\;,
\end{align*}
and are going to show that it is in fact a true martingale.
We study only the second term in \eqref{eq mart term}, the first term is treated similarly with $\dZ$ in place of $\eZ$. Applying Burkholder-Davis-Gundy inequality, we compute, using \eqref{eq ass bound deriv H},
\begin{align*}
\esp{ \sup_{t \in [0,T]} \left|\int_0^t \Gamma_t \dY^i_t\dY^j_t\partial_y a^{ij}_t  \eZ_t\ud W_t\right|} &\le C_\Lambda \esp{(\int_0^T \Gamma_t^2 |\dY_t|^2|\eZ_t|^2 \ud t)^\frac12} 
\\
& \le C_{\Lambda,p} \esp{|\Gamma_T|^p}^\frac1p 
\esp{ \sup_{t\in[0,T]}|\dY_t|^q\left(\int_0^T|\eZ_t|^2 \ud t\right)^q}^\frac1q \,
\\
&\le C_{\Lambda,p}\esp{|\Gamma_T|^p}^\frac1p \,,
\end{align*}
where we used Cauchy-Schwarz inequality, the energy inequality, with the fact that $\sup_t|\dY_t|$ is bounded in any $\LP{r}(\cF_T)$ and $\eZ \in \BP{2}$.
From step 1. we deduce then that the supremum of the local martingale term is integrable and it is thus a martingale.
\\
3. Combining the  results from steps 2.a - 2.c, we get, for $ \mathfrak{B}(\Lambda)$ and $\mathfrak{A}(\Lambda)$ large enough,
%using \eqref{eq useful upper bound},
%\begin{align*}
% \ud [ \Gamma_t \delta Y_t^\dagger A(t,Y_t) \delta Y_t ] \ge&
%\Gamma_t\left(\frac{\epsilon} 6 |\dZ_t|^2 
%-  C_{\Lambda,\epsilon}|(\ef-\zf)(t,\eY_t,\eZ_t)|^2\right) \ud t
%\\
%&- C_\Lambda \Gamma_t \left( |\mathfrak{P}(\zY_t) - \zY_t||\ePhi_t| 
%+  |\mathfrak{P}(\eY_t) - \eY_t||\zPhi_t| \right)
%\ud t + \ud \cM_t\;.
%\end{align*}
\begin{align*}
 \Gamma_t \delta Y_t \cdot A(t,Y_t) \delta Y_t &+ \EFp{t}{\int_t^T \Gamma_s\frac{1}{6L} |\dZ_s|^2 \ud s + \cM_t-\cM_T}
 \\ 
 &\le
 \\
 \EFp{t}{ \Gamma_T \dxi \right. &\cdot \left. A(T,Y_T) \dxi} +
 \\
  C_{\Lambda,K,\sigma} \mathbb{E}_t\left[ \int_t^T \Gamma_s \left\{ |(\ef-\zf)(s,\eY_s,\eZ_s)|^2  \right. \right.&+\left.( |\mathfrak{P}(\zY_s) - \zY_s| +  |\mathfrak{P}(\eY_s) - \eY_s| )(|\ePhi_s| + |\zPhi_s| )
\right\}\ud s \bigg] .
%\\
%+ C_\Lambda \EFp{t}{\int_t^T \Gamma_s (
% |\mathfrak{P}(\zY_s) - \zY_s|\right.&+\left.  |\mathfrak{P}(\eY_t) - \eY_t| )(|\ePhi_t| +
%|\zPhi_t| )
%\ud t } \;.
\end{align*}
Step 2.d allows us to claim that $\EFp{t}{\cM_t-\cM_T}=0$ in the previous inequality.
%where $\cM$ is a local martingale term given by
%\begin{align}
 %\cM_t := \sum_{1\le i,j \le d}  \cM^{ij}_t \; \text{ with } \; \cM^{ij}_t := \int_0^t \Gamma_s \ud M^{ij}_s\;.
%\end{align}
Moreover,  we have, recalling \eqref{eq bound matrix H} and $\Gamma \ge 1$,
\begin{align*}
 |\dY_t|^2 \le L \Gamma_t \delta Y_t \cdot A(t,Y_t) \delta Y_t  \;\text{ and }\; \EFp{t}{ \Gamma_T \dxi \right. &\cdot \left. A(T,Y_T) \dxi} \le L\EFp{t}{ \Gamma_T |\dxi|^2},
\end{align*}
which combined with the previous inequality concludes the proof of \eqref{eq pr stab}.
%\end{bluetext}
%
%
%3. We now show that $\cM$ is a martingale.
%
%Applying Burkholder-Davis-Gundy inequality and observing that $\dY$ is bounded, we compute, using \eqref{eq ass bound deriv H},
%\begin{align*}
%\esp{ \sup_{t \in [0,T]} \left|\int_0^t \Gamma_t \dY^i_t\dY^j_t\partial_y a^{ij}_t  \eZ_t\ud W_t\right|} &\le C \esp{(\int_0^T \Gamma_t^2 |\dY_t|^2|\eZ_t|^2 \ud t)^\frac12} 
%\\
%& \le C \esp{|\Gamma_T|^p}^\frac1p 
%\esp{ \sup_{t\in[0,T]}|\dY_t|^q\left(\int_0^T|\eZ_t|^2 \ud t\right)^q}^\frac1q \,
%\\
%&\le C_p\esp{|\Gamma_T|^p}^\frac1p 
%\end{align*}
%where we used Cauchy-Schwarz inequality, the energy inequality, with the fact that $\sup_t|\dY_t|$ is bounded in any $\LP{r}$ and $\eZ \in \BP{2}$.
%From step 1. we deduce then that the supremum of the local martingale term is integrable, which conclude the  proof for this term. The other terms in \eqref{eq mart term} are treated similarly.
\eproof
%\end{bluetext}

\begin{Remark}
%\begin{enumerate}[i)]
%\item After careful inspection of the previous proof, one could set e.g.
%\begin{align*}
%\mathfrak{A}(\Lambda) :=
%\text{ and }
%\mathfrak{B}(\Lambda) :=
%\end{align*}
%We will not try to characterise optimal constant here, though this can be useful for some applications.
%\item REPHRASE: 
{The dependence upon $\Lambda$ is a key fact that will restrain us to extend straightforwardly to rougher coefficients our main existence and uniqueness results in the non-Markovian case, recall assumption \textbf{(SB)}. %In particular, it does not seem that there is a straightforward generalisation of previous estimates and stability results for an unsmooth function $H$. 
This is a quite important limitation: indeed, when we consider classical applications of obliquely reflected BSDEs to optimal switching problems we have to deal with convex polytopes domains $\cD$ which are corner domains. Thus, in this context $H$ is not smooth enough to apply our results. It justifies Section \ref{se markov} where the Markovian framework is studied under a weaker regularity assumption on $H$ thanks to a different approach.}  
%\end{enumerate}
\end{Remark}

\subsection{Some interesting facts about the class $\mathfrak{T}_\beta$}
\label{subse term xi} 

We first make the following observation.
\begin{Proposition}
\label{prop:distance Zxi Hinfini}
Let $\xi \in \LP{2}(\cF_T)$ satisfying \textbf{(SB)}(i). If we have,  for some $\beta>0$,
\begin{align}\label{eq cond suff}
{d_{\BP{2}}(\cZ^{\xi},\HP{\infty})<\frac{1}{\sqrt{\beta}}}\,,
\end{align}
then $\xi \in \mathfrak{T}_{\beta}$. %\textcolor{red}{reprendre notation}
\end{Proposition}
%\begin{bluetext}
\proof %Let $\tilde{V}$ be in the closure of $\HP{\infty}$ for the $\BP{2}$ norm such that
%$d_{\BP{2}}(\cZ^{\xi},\HP{\infty}) = \NB{2}{\cZ^{\xi}-\tilde{V} }< \frac1{\sqrt{\beta}}$.
We can find $V\in\HP{\infty}$, s.t. $\NB{2}{\cZ^{\xi}-V} = \frac1{(1+\eta)\sqrt{\beta}}$, for some $\eta>0$ small enough. We now set $\lambda := \left(\frac{1+\frac{\eta}{2}}{1+\frac{\eta}{3}}\right)^2\beta$ and we compute, using Young's inequality,
\begin{align*}
|\cZ^\xi|^2\le(1+\frac{\eta}3)|\cZ^\xi-V|^2 + (1+\frac3{\eta})|V|^2.\
\end{align*}
This leads, using Hölder inequality, to
\begin{align*}
\esp{e^{\lambda\int_0^T|\cZ^\xi_t|^2 \ud t}}
\le
C\esp{e^{(1+\frac{\eta}3)^2\lambda\int_0^T|\cZ^\xi_t-V_t|^2\ud t}}^{\frac1{1+\frac{\eta}3}}
\end{align*}
where we used the fact that $V \in \HP{\infty}$. Since $\NB{2}{(1+\frac{\eta}3)\sqrt{\lambda}(\cZ^\xi-V)} = \frac{1+\frac{\eta}2}{1+\eta} < 1$, we can apply the John-Nirenberg inequality, see Theorem 2.2 in \cite{Kazamaki-94}, to obtain
\begin{align*}
\esp{e^{(1+\frac{\eta}3)^2\lambda\int_0^T|\cZ^\xi_t-V_t|^2\ud t}}<\infty\;,
\end{align*}
which concludes the proof.
\eproof

Proposition \ref{prop:distance Zxi Hinfini} only suggests a sufficient condition. In the case $\beta = +\infty$, for which condition \eqref{eq cond suff} should read $d_{\BP{2}}(\cZ^{\xi},\HP{\infty})=0$, it is known that the condition is not necessary. We refer the interested reader to the paper \cite{Schachermayer-96}, where this question is treated with more details.

\vspace{5pt}

The next result shows that a class of path-dependent function of some smooth processes are naturally contained in $\mathfrak{T}_\beta$ and actually for all $\beta > 0$. This class is quite important for applications. 
%This Corollary follows mainly from the next Proposition.

\begin{Proposition}
\label{prop:g(X) dans Hinfini}
Let $X \in \SP{2}$ such that for all $t,s \le T$, the Malliavin derivatives of $X_s$ denoted $D_tX_s$ is well defined and satisfies $\NL{\infty}{\sup_{t,s} |\EFp{t}{D_t X_s}|}< \infty$. Let $g:\cC^0([0,T],\R^n) \rightarrow \R^d $ be a uniformly continuous function, then denoting $\xi = g\left((X_s)_{s\in[0,T]}\right)$, we have that $\cZ^\xi \in \widebar{\HP{\infty}}^{\BP{2}}$.
\end{Proposition}

%\begin{bluetext}
1.a %\textcolor{magenta}{reprendre pour virer la bornitude} 
We first start by considering a sequence $(g_N)$ of $N$-Lipschitz regularisation of $g=(g^1,\dots,g^d)$ given by
$$g^i_N (x) = \inf_{u \in \cC^0([0,T],\R^n)} \{ g^i(u)+N \|u-x\|_\infty \}\;,\quad \text{ for all } x \in \cC^0([0,T],\R^n), \quad 1 \le i \le n.$$
Let us observe that $g_N$ is finite for $N$ large enough due to the linear growth of $g$.  Then 
we have, for all $ x \in \cC^0([0,T],\R^n)$ and $1 \le i \le n$,
\begin{align*}
 g^i(x) \ge g_N^i(x) &\ge \inf_{u \in \cC^0([0,T],\R^n)} \{ g^i(x)-\omega_{g^i}(|u-x|_\infty)+N |u-x|_\infty \}\\
 &\ge g^i(x) + \inf_{u \in \cC^0([0,T],\R^n)} \{ N|u|_\infty-\omega_{g^i}(|u|_\infty) \}
\end{align*}
where $\omega_{g^i}$ is a concave modulus of continuity for the uniformly continuous component $g^i$ of $g$. Thus we get
\begin{align}
\label{eq conv approx unif}
 |g_N-g|_{\infty} \le C\sum_{i=1}^d \sup_{u \in \cC^0([0,T],\R^n)} \{ \omega_{g^i}(|u|_\infty)-N|u|_\infty \}:=c(N).
\end{align}
Since $\omega_{g^i}(h)=o(1)$ when $h \rightarrow 0^+$ then $c(N) = o(1)$ when $N \rightarrow + \infty$.
% We first show that we there exists $g_N$ a $N$-Lipschitz regularisation of $g=(g^1,\dots,g^d)$ satisfying
% \begin{align}\label{eq conv approx unif}
% \sup_{x \in \cC^0([0,T],\R^n)} |g(x)-g_N(x)| \le C\sum_{i=1}^d\frac1N+\omega_{g^i}\!\left(\frac{2|g|_\infty}N\right) =:c(N)
% \end{align}
% where $\omega_{g^i}$ is a modulus of continuity for the U.C. component $g^i$ of $g$.
% We work component by component and w.l.o.g we consider only the first one for the rest of the proof. Let us define
% $$g^1_N (x) = \inf_{u \in \cC^0([0,T],\R^n)} \{ g^1(u)+N |u-x|_\infty \}\;,\quad \text{ for all } x \in \cC^0([0,T],\R^n).$$
% We observe then, that for all $N\ge1$, there exists $x_N$ such that
% \begin{align}\label{eq temp N lip reg}
% g^1_N(x)\le g^1(x_N) + N|x_N-x|_\infty \le g^1(x)+\frac1N 
% \end{align}
% and
% \begin{align}\label{eq temp N lip reg 2}
%  g^1_N(x) + \frac1N \ge g^1(x_N) + N|x_N-x|_\infty.
% \end{align}
% Using \eqref{eq temp N lip reg}, we obtain, for $N$ large enough, that
% \begin{align*}
% |x_N-x|_\infty \le \frac{3|g|_\infty}N \;\text{ and } -\frac1N \le g^1(x)-g^1_N(x).
% \end{align*}
% Equation \eqref{eq temp N lip reg 2} yields that
% \begin{align*}
% g^1_N(x) - g^1(x) \ge  g^1(x_N) - g^1(x) - \frac1N \;.
% \end{align*}
% Combining the three previous inequalities, we get
% \begin{align*}
% |g^1_N(x) - g^1(x)| \le \frac1N + \omega_{g^1}\!\!\left(\frac{3|g|_\infty}N\right)\;.
% \end{align*}

1.b Defining $\cY^N_t := \EFp{t}{g_N(X_\cdot)} = g_N(X_\cdot) - \int_t^T \cZ^N_s \ud W_s$ and applying It\^o's formula to $|\cY^N_\cdot-\cY^\xi_\cdot|^2$, we compute
\begin{align*}
|Y^N_t-\cY^\xi_t|^2 + \EFp{t}{\int_t^T|\cZ^N_t-\cZ^\xi_t|^2 \ud t} = \EFp{t}{|g_N(X_\cdot)-g(X_\cdot)|^2} \le c(N)^2
\end{align*}
recall \eqref{eq conv approx unif}.
From this, we deduce that for all $\epsilon >0$, there exists $N_\epsilon$, s.t. for all $N \ge N_\epsilon$,
\begin{align}\label{eq reg first step}
\NB{2}{\cZ^N-\cZ^\xi} \le \epsilon\;.
\end{align}
2. We now show that $\cZ^N$  introduced above, belongs to  $\HP{\infty}$. This fact combined with \eqref{eq reg first step} proves the statement of the proposition.\\
Following Lemma 4.1 in \cite{Ma-Zhang-02} there exists a family $\Pi = \{\pi\}$ of partitions of $[0,T]$ and a family of discrete functionals $\{g_{N,\pi}\}$ such that 
\begin{itemize}
 \item for each $\pi \in \Pi$, with $\pi : 0=t_0 <...<t_m=T$, we have that $g_{N,\pi}\in C^{\infty}_b(\mathbb{R}^{d(m+1)})$, and satisfies
 \begin{align}\label{eq bound deriv}
 \sum_{i=0}^m |\partial_{x_i} g_{N,\pi}(x)| \leqslant N, \quad \forall x \in  \cC^0([0,T],\R^n),
 \end{align}
 where $g_{N,\pi}(x):=g_{N,\pi}(x(t_0),...,x(t_m))$.
 \item for any $x \in  \cC^0([0,T],\R^n)$ it holds that
 \begin{align}\label{eq approx disc}
 \lim_{|\pi| \rightarrow 0} |g_{N,\pi}(x)-g_N(x)| = 0.
 \end{align}
\end{itemize}
We naturally consider $(\cY^{N,\pi},\cZ^{N,\pi})$ given by
\begin{align*}
\cY^{N,\pi}_t := \EFp{t}{g_{N,\pi}(X) } = g_{N,\pi}(X) - \int_t^T\cZ^{N,\pi}_s \ud W_s, \; \quad 0 \leqslant t \leqslant T.
\end{align*}
2.a By the Clark-Ocone formula, we have that 
\begin{align*}
\cZ^{N,\pi}_t &= \EFp{t}{D_tg_{N,\pi}(X)} 
\\
&= \sum_{i=1}^m \partial_{x_i}g_{N,\pi}(X) \EFp{t}{D_t X_{t_i} }.
\end{align*}
Now, using \eqref{eq bound deriv} and the assumption on $(D_tX_s)_{t,s \le T}$, we obtain
\begin{align}\label{eq sinfinity bound on ZNpi}
\NS{\infty}{\cZ^{N,\pi}}\le C_N\;.
\end{align}
\\
2.b {Combining \eqref{eq approx disc} and the fact that $|g_N(X) - g_{N,\pi}(X)| \leqslant C_N(1+\sup_{t \in [0,T]} |X_t|) \in L^2$, we can use the dominated convergence theorem %\textcolor{red}{ je pense que $g$ a croissance lineaire et ok pour cela et $X$ dans $\SP{2}$...}, we 
to get}
\begin{align*}
 \lim_{|\pi| \rightarrow 0} \esp{|g_N(X)-g_{N,\pi}(X)|^2} = 0\;,
\end{align*}
which leads to
\begin{align*}
\lim_{|\pi| \rightarrow 0}\esp{\int_0^T|\cZ^N-\cZ^{N,\pi}|^2\ud t} = 0\;.
\end{align*}
Up to a subsequence, we have
$\cZ^{N,\pi}\rightarrow \cZ^N$ $\ud \P \otimes \ud t$-a.e. and, moreover, $|\cZ^{N,\pi}|\le C_N$, recall \eqref{eq sinfinity bound on ZNpi}.  We thus obtain for (a version of) the limit process 
\begin{align}
\int_0^T|\cZ^{N}_t|^2\ud t \le TC_N^2 \;, \quad \P-a.s.
\end{align}
which concludes the proof of this step.

3. Finally, we remark that $\cZ^\xi \in \BP{2}$ since
$\cZ^\xi -\cZ^N \in \BP{2}$ and $\cZ^N \in \HP{\infty} \subset \BP{2}$. We conclude the proof by using \eqref{eq reg first step}.
\eproof
%\end{bluetext}

\vspace{5pt}
{We obtain the following direct corollary, which gives a sufficient condition on models in a path-dependent framework to check the admissibility of the terminal condition.}

\begin{Corollary}\label{co useful app}
Let $X$ be solution of the Lipschitz SDE 
$$X_t = x+\int_0^t b(X_s)ds+\int_0^t \sigma(X_s) dW_s,$$
where $\sigma$ and $b$ are Lipschitz continuous functions and $\sigma$ is bounded.  \\
Set $\xi := g((X_s)_{t \in [0,T]})$ where $g$ is a uniformly continuous function on $\cC([0,T],\R^d)$, then $\xi$ belongs to $\mathfrak{T}_\beta$, for all $\beta > 0$. Moreover, if $\tilde{\xi} \in \LP{\infty}$, then $\xi + \tilde{\xi}$ belongs to $\mathfrak{T}_{\beta}$ for all {$\beta < \NL{\infty}{\tilde{\xi}}^{-2}$}.
\end{Corollary}
%\begin{bluetext}

{\proof  
When $\sigma$ and $b$ are smooth enough, it is well known, see e.g. \cite{Nualart-06}, that $X$ is Malliavin differentiable and, for all $1 \le i \le k$, $(D^i_tX_s)_{s \in [t,T]}$ is solution of the linear SDE given by 
$$D^i_t X_s = \sigma^i(X_t) +\int_t^s \nabla b(X_r) D^i_t X_r \ud r + \int_t^s \sum_{j=1}^k \nabla \sigma^j(X_r)  D^i_t X_r \ud W^j_r,\quad t \leqslant s \leqslant T .$$
Then, we easily get that $|\mathbb{E}_t[D_t X_s]| \le e^{K_b T} M$ with $K_b$ the Lipschitz constant of $b$ and $M$ a bound of $\sigma$. Thus we can apply Proposition \ref{prop:g(X) dans Hinfini} to get the first part of the result. When coefficients are not smooth enough, a standard approximation gives us the result, pointing out the fact that $\NL{\infty}{\sup_{t,s} |\EFp{t}{D_t X_s}|}$ can be uniformly bounded with respect to the approximation.
For the second part of the corollary, we just have to remark that 
$$d_{\BP{2}}(\cZ^{\xi+\tilde{\xi}},\HP{\infty})\le \NB{2}{\cZ^{\xi+\tilde{\xi}}-\cZ^{\xi}} + d_{\BP{2}}(\cZ^{\xi},\HP{\infty}) = \NB{2}{\cZ^{\tilde{ \xi}}} .$$
Moreover, applying It\^o's formula to $|\cY^{\tilde{\xi}}_t|^2$, we compute
$$\NB{2}{\cZ^{\tilde{\xi}}} \le \NL{\infty}{\tilde{\xi}},$$
which implies
$$d_{\BP{2}}(\cZ^{\xi+\tilde{\xi}},\HP{\infty})\le  \NL{\infty}{\tilde{\xi}}.$$
Thus, we just have to apply Proposition \ref{prop:distance Zxi Hinfini} to conclude.
\eproof}
%\end{bluetext}
%}

%!TEX root = main.tex
\section{Existence and uniqueness in a regular setting}
\label{se smooth setting}

In this section, we obtain an existence and uniqueness result in  a non Markovian  setting, working under assumption \textbf{(SB)} and  considering terminal condition in the class $\mathfrak{T}_\beta$, for some $\beta>0$. This $\beta$, as shown in the previous section depends dramatically on the smoothness of the coefficients. Our proof is done in two  main steps. In the first step, we restrict to the case of a bounded terminal condition. We study the wellposedness of the penalised equation, and prove their convergence to an obliquely reflected BSDE. In a second step, we extend our result to all terminal condition in the class $\mathfrak{T}_\beta$.
%, for some $\beta>0$ depending on the smoothness of the coefficients.
%a smooth domain and bounded coefficient framework. %Precisely, we shall work in the setting of assumption \textbf{(SB)}.

%\textcolor{green}{
%\begin{Remark} In our setting of bounded coefficients, the extension of our result to unbounded domain is straightforward: It suffices to consider a smooth bounded domain $\tilde{\cD}$ s.t. $ \partial \cD \cap B(0,2R) \subset \partial \tilde{\cD}$, where $R$ is large enough according to Proposition \ref{estimee a priori 1}, equation \eqref{eq estimate penalised bounded setting} and Proposition \ref{pr a priori orbsde}. On the frontier $\cF:=\partial \tilde{\cD} \setminus \partial \cD$, we can extend smoothly the function $H$ sucht that $\textbf{(SB)}(iii)$  is satisfied (by considering first that $H=I_d$ and then using a regularisation argument). Note that any extension is possible since only $\partial \cD \cap B(0,R)$ can be visited by the solution $Y$.
%%\textcolor{red}{extension of the direction of reflection}
%\end{Remark}
%%\textcolor{red}
%{
%remarque: il semble possible d'obtenir un resultat plus general en terme de bornitude...
%il faut verifier que la propriete d'integrabilite du $\cZ^\xi$ passe a la solution de l'edsr reflechie.
%avec une propriete de moment expo du y, cela permettrait d'obtenir les estimees expo sur le Phi (cf proposition bornitude RBSDE partie 1 preuve etc.
%}
%}

%\subsection{Existence and uniqueness of the penalized BSDE}
\label{section:exist-uniq penalized BSDE}

% To obtain an existence and uniqueness result for the penalized BSDE \eqref{eq:edsr penalized gene}, we assume some regularity assumptions on $f$ and $H$.
% \paragraph{Assumption (ALip)}
%$ $
% There exist constants $K_y$ and $K_z$ such that, for all $(t,y,y',z,z') \in [0,T] \times \mathbb{R}^d \times\mathbb{R}^d \times \mathbb{R}^{d \times k}\times \mathbb{R}^{d \times k}$ we have
% $$|\bar{H}(t,y,z)-\bar{H}(t,y',z')|+|f(t,y,z)-f(t,y',z')| \leqslant K_y|y-y'|+K_z|z-z'|.$$
% 

\subsection{Bounded terminal condition}
We first obtain some results on the penalised BSDE that will be used later in this section and also in Section \ref{se markov} in the Markovian case. We thus essentially work here under the assumption \textbf{(A)}.

\vspace{5pt}
\noindent We start with the following lemma that  verifies the well-posedness of equation \eqref{eq:edsr penalized gene} under some classical conditions.% is well posed.

\begin{Lemma}
\label{le:euEDSRpenalise} %\textcolor{red}{a deplacer?}
We assume that \textbf{(A)} is in force and that $f$ and $H$ are Lipschitz continuous with respect to $(y,z)$. Then, for all $n \in \mathbb{N}$ there exists a unique solution to \eqref{eq:edsr penalized gene} in $\mathscr{S}^2\times\mathscr{H}^2$. 
\end{Lemma}
%\textcolor{red}{erreur? manque un facteur $2$ dans la derivee... continuite??}

\proof
% We remark that 
% $$|f(t,x,y,z)| \leq C(1+|x|^p + |y| + |z|) \quad \textrm{and} \quad |H(t,x,y,z)\nabla \varphi_n(y)| \leq C_n(1+|y|).$$
% Moreover, for all $(t,x) \in [0,T] \times \R^q$, $(y,z) \mapsto f(t,x,y,z)+H(t,x,y,z)\nabla \varphi_n(y)$ is a continuous fonction: indeed $H$ is a bounded function, continous on $\R^d \setminus \bar{\cD}$ and $\nabla \varphi_n(y)=0$ for $y\in \bar{\cD}$. So a direct application of Theorem 27.2 in \cite{Hamadene-Lepeltier-Peng-97} gives us the result (see also remark 27.3 ii)).
Since $\mathcal{D}$ is convex, $\varphi^M_n$ is convex and $nM$-Lipschitz continuous, recall \eqref{de varphi n}. Indeed,
denoting $\mathcal{D}_M := \set{ y \in \R^d | d(y,\mathcal{D}) \le M}$, we have that 
% \begin{align*}
% \forall y \in \mathcal{D}_M, \varphi_M(y) = d^2(y,\mathcal{D}) \;.
% \end{align*}
 \begin{align}\label{eq expression  varphi}
 \varphi^M_n(h) =
\left \{
\begin{array}{lcl}
 n  \frac12d^2(h,\mathcal{D}) & \text{if} & h \in \mathcal{D}_M
 \\
n Md(h,\mathcal{D})  -\frac{nM^2}2& \text{if}  & h \notin \mathcal{D}_M
\end{array}
\right.
\end{align}
and
\begin{align}\label{eq expression nabla varphi}
  \nabla \varphi^M_n(y) =
\left \{
\begin{array}{lcl}
0 & \text{if} & y \in \bar{\cD}\\
 n d(y,\cD)\frac{y-\mathfrak{P}(y)}{|y-\mathfrak{P}(y)|} & \text{if} & y \in \mathcal{D}_M \setminus \bar{\cD}
 \\
nM \frac{y-\mathfrak{P}(y)}{|y-\mathfrak{P}(y)|} & \text{if}  & y \notin \mathcal{D}_M
\end{array}
\right. .
\end{align}
Finally ${H}$ and $\nabla \varphi^M_n$ are two Lipschitz 
bounded functions which proves that the penalised BSDE
\eqref{eq:edsr penalized gene} has a Lipschitz driver: the classical result of 
\cite{Pardoux-Peng-90} then applies to get the existence and uniqueness result.
\eproof

\begin{Lemma} \label{le control Phi and dist penalised} Assume that \textbf{(A)} holds and that there exists a solution to \eqref{eq:edsr penalized gene} in $\SP{2}\times\HP{2}$. Then, $(Y^n,Z^n,\nabla \varphi_n^M(Y^n))$ 
satisfies Condition \eqref{eq hyp gen Phi} for some $K>0$ and  for some $c>0$ we have
\begin{align}\label{eq lemma control Phi penalised}
\sup_{t \in [0,T]} \esp{  \varphi^M_n(Y_t^{n})}+\esp{\int_0^T |\nabla \varphi^M_n (Y_s^n)|^2 \ud s } \leqslant c \mathbb{E}\left[|\xi|^2+\int_0^T |\alpha_s|^2 ds\right].
\end{align}
Importantly, $K$ and $c$ do not depend on $n$, {nor $M$}.\\ 
Moreover, if \textbf{(SB)} holds, then, there exists $c':=c'(\sigma^\xi)$, which does not depend on $n$ {nor $M$}, such that
\begin{align}\label{eq control Phi penalised under sb}
\sup_{t\in[0,T]} \varphi^M_n(Y_t^{n}) + \NB{2}{\nabla \varphi_n^M(Y^n)}^2 \le c'\;.
\end{align}
\end{Lemma}

%\begin{bluetext}
\proof Since $\varphi^M_n$ is a $C^1$ convex function, we have  the following inequality (see Lemma 2.38 in \cite{pardoux2014stochastic}): for $s \in [t,T]$, 
\begin{align}
\label{Ito fonction convexe nm}
 \varphi^M_n(Y^{n}_s) 
& +  \int_s^T  \nabla \varphi^M_n(Y^{n}_u) \cdot {H}(u,Y^{n}_u,Z_u^{n})  \nabla \varphi^M_n(Y^{n}_u)  \ud u\\
\le &  \varphi^M_n(\xi) +  \int_s^T  \nabla \varphi^M_n(Y^{n}_u) \cdot f(u,Y^{n}_u,Z^{n}_u)  \ud u
-  \int_s^T  \nabla \varphi^M_n(Y^{n}_u) \cdot Z^{n}_u \ud W_u , \nonumber
\end{align}
and we recall that $\varphi^M_n(\xi)=0$. We observe, using \eqref{inegalite H a intro} that
\begin{equation}
   \nabla \varphi^M_n(Y^{n}_u) \cdot {H}(u,Y^{n}_u,Z^{n}_u)  \nabla \varphi^M_n(Y^{n}_u) 
\ge
\eta  |\nabla \varphi^M_n(Y^{n}_u)|^2
\label{inegalite 1  nm}
\end{equation}
and   combining Cauchy-Schwarz inequality with Young's inequality
\begin{align*}
 \int_s^T  \nabla \varphi^M_n(Y^{n}_u) \cdot f(u,Y^{n}_u,Z^{n}_u)  \ud u
 \le \frac{\eta}2  \int_s^T |\nabla \varphi^M_n(Y^{n}_u)|^2 \ud u + \frac2\eta \int_s^T |f(u,Y^{n}_u,Z^{n}_u) |^2 \ud u\;.
\end{align*}
From this, we deduce
\begin{align}\label{eq control penalised temp 1}
 \varphi^M_n(Y^{n}_t)  + \EFp{t}{\int_t^T |\nabla \varphi^M_n(Y^{n}_u)|^2 \ud u}\le \frac4{\eta} \EFp{t}{\int_t^T |f(u,Y^{n}_u,Z^{n}_u) |^2 \ud u}\;, 
\end{align}
which proves \eqref{eq hyp gen Phi} for $(Y^n,Z^n,\Phi^n)$.
This allows then to invoke Lemma \ref{le first basic a priori} to obtain \eqref{eq lemma control Phi penalised} under \textbf{(A)}.
Under \textbf{(SB)}, \eqref{eq control penalised temp 1} allows also to conclude recalling that $f$ is Lipschitz continuous, $\theta^\xi \in \BP{2}$ and \eqref{eq control DY DZ}.
\eproof
%\end{bluetext}

\vspace{5pt}
%\subsection{Existence and uniqueness for the obliquely reflected BSDE}

We now prove our first existence result for the obliquely reflected BSDE 
\begin{equation}
 \label{eq smooth and bounded ORBSDE}
 \left\{ \begin{aligned} &Y_t = \xi+\int_t^T f(s,Y_s,Z_s)\ud s-\int_t^T H(s,Y_s)\Phi_s \ud s-\int_t^T Z_s \ud W_s, \quad 0 \leqslant t \leqslant T,\\
         & Y_t \in \bar{\cD}, \quad \Phi_t \in\partial \varphi(Y_t),\quad 
          \int_0^T  \1_{\set{Y_t \notin \partial \cD}}  |\Phi_t| \ud t= 0.
        \end{aligned}
 \right.
\end{equation}
%in this smooth and bounded setting and, in particular, under the assumption $\HYP{\xi}$.

\begin{Proposition} \label{pr smooth and bounded exist orbsde} Assume that \textbf{(SB)} holds and that $\xi \in \LP{\infty}\cap\mathfrak{T}_{\mathfrak{B}(\Lambda)}$. Then, there exists a  solution in  $\SP{2}\times\HP{2}\times\HP{2}$ to the obliquely reflected BSDE \eqref{eq smooth and bounded ORBSDE}. 
%Moreover, it satisfies
%\begin{align*} 
%  \NS{\infty}{Y}  + 
 % \NB{2}{Z}
 % +
 % \NB{2}{\Phi}
  %\le 
 % C(\NL{\infty}{\xi} + \NB{2}{\alpha}) < \infty \, ,
%\end{align*}
%and
%\textcolor{red}{condition integrabilite exp pour Z => vraie par Lemme 2.1 en fait...}
\end{Proposition}

%\begin{bluetext}
\proof
To obtain the existence result, we consider a sequence of penalised BSDEs given by equation \eqref{eq:edsr penalized gene} for which we have existence and uniqueness from Lemma \ref{le:euEDSRpenalise}. In the definition of $\varphi^M_n$, recall \eqref{de varphi n}, we set 
$M = 2 c$ where $c$ is given in Corollary \ref{co Y bounded}. In particular, we observe that for this choice of $M$, for $0 \le t \le T$,
\begin{align}\label{eq shape of Phi}
\Phi^n_t := \nabla \varphi^M_n(Y^n_t) = n \left(Y^n_t-\mathfrak{P}(Y^n_t) \right)
\text{ and } \frac1n\varphi_n(Y^n_t) = {\frac12}|\mathfrak{P}(Y^n_t) - Y^n_t|^2  
\;,
\end{align}
recall \eqref{eq expression varphi} and \eqref{eq expression nabla varphi}. We will use this fact later on. 
\\
1.a. We now prove that $(Y^{n},Z^{n})$ is a Cauchy sequence in $\SP{2}\times\HP{2}$. Indeed, let $m \ge 0$ and $n \ge 0$, thanks to Lemma \ref{le control Phi and dist penalised} we can apply Proposition \ref{pr stability} to obtain
\begin{align}\label{eq from prop 2.3}
 &\sup_{t \in [0,T]} \esp{|Y^n_t-Y^m_t|^2} + \NH{2}{Z^n-Z^m}^2
 \\ &\le C_{\Lambda} 
\esp{ \int_0^T \Gamma^{n,m}_s \left( |\mathfrak{P}(Y^n_s) - Y^n_s| + |\mathfrak{P}(Y^m_s) - Y^m_s|\right)\left(|\Phi^m_s| 
 +  |\Phi^n_s| \right)
\ud s\;.
} =: A^{n,m}. \nonumber
\end{align}
Let us notice that, from Proposition \ref{pr stability} again, there exist $p>1$ and a constant $C$ such that
\begin{align}\label{eq control Gamma unif}
\esp{|\Gamma^{n,m}_T|^p} \le C\;,
\end{align}
where, importantly, $p$ and $C$ do not depend on $(n,m)$.
Applying It\^o's formula to $|Y^n-Y^m|^2$ on $[0,T]$, we compute, using usual arguments,
\begin{align*}
\NS{2}{Y^n-Y^m}^2 &\le C\esp{\int_0^T|Y^n_t-Y^m_t|\left(|\Phi^n_t|+|\Phi^m_t| \right)\ud t}
\\
&+C\esp{\sup_{t \in [0,T]}\left|\int_0^t (Y^n_s-Y^m_s)(Z^n_s-Z^m_s)\ud W_s\right|} \;.
\end{align*}
Using Burkholder-Davis-Gundy inequality and Young's inequality, we obtain
\begin{align*}
\NS{2}{Y^n-Y^m}^2 \le C\esp{\int_0^T|Y^n_t-Y^m_t|\left(|\Phi^n|+|\Phi^m_t| \right)\ud t}
+C\NH{2}{Z^n-Z^m}^2 \;.
\end{align*}
Applying Cauchy-Schwarz inequality, and using Lemma \ref{le control Phi and dist penalised}, we  get
\begin{align}
\NS{2}{Y^n-Y^m}^2  \le C \left( \NH{2}{Y^n-Y^m} +  \NH{2}{Z^n-Z^m}^2 \right)\;.
\end{align}
Combining the previous inequality with \eqref{eq from prop 2.3}, we  have 
\begin{align}\label{eq end step 2.a}
\NS{2}{Y^n-Y^m}^2 + \NH{2}{Z^n-Z^m}^2 \le C\left( A_{n,m} + \sqrt{A_{n,m} } \right).
\end{align}
%\textcolor{red}{bien verifier l'independance de la constante en $m$ et $n$}
1.b We now study the $A_{n,m}$ term. 
%In the definition of the $(\varphi_n)$, see equation \ref{de varphi n} and also Remark \ref{re quand il y a choix il y a bon choix} , we choose 
%\begin{align}\label{eq choice of M}
%M := 2C(\NL{\infty}{\xi} +  \NB{2}{\alpha}), 
%\end{align}
%recall Lemma \ref{le bmo setting estimate}. With this choice, we observe that 
%\begin{align*}
%{\frac12}|\mathfrak{P}(Y^n_s) - Y^n_s|^2 = \frac1n\varphi_n(Y^n_s) \quad \text{ and } %\quad
%{\frac12} |\mathfrak{P}(Y^m_s) - Y^m_s|^2 = \frac1m\varphi_m(Y^m_s)\;.
%\end{align*}
We first observe, recalling Lemma \ref{le control Phi and dist penalised} and \eqref{eq shape of Phi},
%\eqref{eq control Phi penalised under sb},
\begin{align*}
\esp{ \int_0^T \Gamma^{n,m}_t  |\mathfrak{P}(Y^n_t) - Y^n_t||\Phi^m_s| \ud s}
&\le \NL{\infty}{\sup_{t} |\mathfrak{P}(Y^n_t) - Y^n_t|} \esp{\Gamma^{n,m}_T \int_0^T
|\Phi^m_s | \ud s}
\\
&\le \frac{C}{\sqrt{n}}\esp{\Gamma^{n,m}_T \int_0^T
|\Phi^m_s | \ud s}\;.
\end{align*}
Applying H\"older inequality, denoting $q$ the conjugate exponent of $p$ introduced in \eqref{eq control Gamma unif}, we deduce from the previous inequality
\begin{align*}
\esp{ \int_0^T \Gamma^{n,m}_t  |\mathfrak{P}(Y^n_t) - Y^n_t||\Phi^m_s| \ud s}
&\le \frac{C}{\sqrt{n}} \esp{\left(\int_0^T |\Phi^m_s|^2\ud s\right)^{\frac{q}2} }.
\end{align*}
%recall \eqref{eq control Gamma unif}.
Then, combining the energy inequality with \eqref{eq control Phi penalised under sb}, we conclude
\begin{align}
\esp{ \int_0^T \Gamma^{n,m}_t  |\mathfrak{P}(Y^n_t) - Y^n_t||\Phi^m_s| \ud s}
&\le \frac{C}{\sqrt{n}} \;.
\end{align}
Similarly we obtain,
\begin{align*}
\esp{ \int_0^T \Gamma^{n,m}_t  |\mathfrak{P}(Y^m_t) - Y^m_t||\Phi^n_s| \ud s}
\le \frac{C}{\sqrt{m}}
\end{align*}
and
\begin{align*}
\esp{ \int_0^T \Gamma^{n,m}_t ( |\mathfrak{P}(Y^n_t) - Y^n_t|+|\mathfrak{P}(Y^m_t) - Y^m_t|)|\Phi^m_s| \ud s}
\le C\left(\frac{1}{\sqrt{m}} + \frac{1}{\sqrt{n}} \right).
\end{align*}

Combining the previous inequalities with \eqref{eq end step 2.a}, we compute that
\begin{align*}
\NS{2}{Y^n-Y^m}^2 + \NH{2}{Z^n-Z^m}^2 \le C\left( n^{-\frac14} + m^{-\frac14} \right)\,,
\end{align*}
which proves that $(Y^n,Z^n)_n$ is a Cauchy sequence in $\SP{2}\times\HP{2}$. We denote $({Y},{Z})$ its limit.
\\
2. We now prove that $(Y,Z)$ is solution to an obliquely reflected BSDE, namely we pass to the limit in \eqref{eq:edsr penalized gene}. Let us first observe that, passing to the limit in \eqref{eq control Phi penalised under sb} yields that $Y \in \bar{\cD}$ as expected.
\\
2.a We first study the reflecting term. 
Since, by Lemma  \ref{le control Phi and dist penalised}, 
$$  \esp{\int_0^T|\nabla \varphi_n(Y^n_s)|^2 \ud s} \le C,$$
we have, up to a subsequence, the following weak $L^2([0,T]\times\Omega)$-convergence:
%$$\nabla \varphi_n(Y^n_\cdot ) \rightharpoonup \Psi, \quad \text{when } n \rightarrow + \infty,$$
%and
$$\nabla \varphi_n(Y^n) \rightharpoonup \Phi, \quad \text{when } n \rightarrow + \infty.$$
%Using Mazur's Lemma, we know that there exists convex combination of the above %converging strongly in $L^2([0,T]\times\Omega)$, namely
%\begin{align*}
% {}^p\!S^{\Psi} := \sum_{r=p}^{N_p}\lambda^p_r \nabla \varphi_r(Y^r_{\cdot})\mathbbm{1}_{Y_\cdot \in \mathcal{O}} \stackrel{p \rightarrow \infty}{\rightarrow} \Psi
%\end{align*}
%and 
%\begin{align*}
% {}^p\!S^{\Phi} := \sum_{n=p}^{N_p}\lambda^p_n \nabla \varphi_n(Y^n_{\cdot}) %\stackrel{p \rightarrow \infty}{\rightarrow} \Phi,
%\end{align*}
%where $\lambda^p_n \ge 0$ for all $p \in \mathbb{N}$ and $p \le n \le N_p$, and $%\sum_{n=p}^{N_p}\lambda^p_n=1$.
%\\
Let $(V_t)_{t \in [0,T]}$ be a continuous adapted process valued in $\bar{\cD}$. From the convexity property of $\cD$ and the fact that $\nabla \varphi_n(Y^n) = n(Y^n-\mathfrak{P}(Y^n_s))$, recall \eqref{eq shape of Phi}, we have
\begin{align*}
\int_0^T(Y^n_t-V_t)^\dagger\nabla \varphi_n(Y^n_t) \ud t \le 0 \;.
\end{align*}
By strong convergence of $(Y^n)_{n \ge 0}$ to $Y$, weak convergence of $(\nabla \varphi_n(Y^n_{\cdot}))_{n \ge 0}$ and the uniform $L^2$-bound on $ \nabla \varphi_n(Y^n_{\cdot})$, recall Lemma  \ref{le control Phi and dist penalised}, we obtain
\begin{align*}
\esp{\int_0^T(Y_t-V_t)^\dagger \Phi_t\ud t \1_A} \le 0\,,
\end{align*}
for all $A \in \cF_T$. This leads to
%\begin{align*}
$ \int_0^T(Y_t-V_t)^\dagger \Phi_t\ud t \le 0\;. $
%\end{align*}
Using Lemma 2.1 in \cite{Gegout-Petit-Pardoux-96} $\omega$-wise, we obtain that
\begin{align*}
\Phi \in \partial \varphi(Y) \; \text{ and } \;  \int_0^T  \1_{\set{Y_t \notin \partial \cD}}  |\Phi_t|\ud t = 0 \;,
\end{align*}
which fully characterise $\Phi$.

2.b Now we want to show that $(Y,Z,\Phi)$ is solution of \eqref{eq smooth and bounded ORBSDE}. By strong convergence of $(Y^n,Z^n)$ to $(Y,Z)$ and the Lipschitz-continuity of $f$, we have
\begin{align*}
f(\cdot,Y^n_\cdot, Z^n_\cdot) \overset{\HP{2}}{\longrightarrow} f(\cdot,Y_\cdot, Z_\cdot)
\;\text{ and } \;\int_0^t Z^n_s \ud W_s \overset{\LP{2}}{\longrightarrow} \int_0^t Z_s \ud W_s \;, 
\end{align*}
for all $t \le T$. Moreover, $\Phi^n \rightharpoonup \Phi$ in $L^2([0,T]\times\Omega)$, when $n \rightarrow + \infty$. Using Mazur's Lemma, we know that there exists a convex combination of the above converging strongly in $L^2([0,T]\times\Omega)$, namely
\begin{align*}
 {}^p\!{\Phi} := \sum_{r=p}^{N_p}\lambda^p_r \Phi^r \stackrel{p \rightarrow \infty}{\rightarrow} \Phi,
\end{align*}
where $\lambda^p_r \ge 0$ for all $p \in \mathbb{N}$ and $p \le r \le N_p$, and $\sum_{r=p}^{N_p}\lambda^p_r =1$.
Let us observe that by strong convergence, the following combination
\begin{align*}
 ({}^p\!Y,{}^p\!Z) := \sum_{r=p}^{N_p}\lambda^p_r (Y^r,Z^r)
\end{align*}
still converges to $(Y,Z)$ in $\SP{2} \times \HP{2}$ and, by strong convergence, 
\begin{align*}
 \sum_{r=p}^{N_p}\lambda^p_r f(\cdot,Y^r,Z^r) &\overset{\HP{2}}{\longrightarrow}
  f(\cdot,Y,Z)\; \text{ and } \; \int_0^t {}^p\!Z_s \ud W_s \overset{\LP{2}}{\longrightarrow} \int_0^t Z_s \ud W_s\;, t \le T \;.
\end{align*}
Moreover, we remark that
 \begin{align*}
 \sum_{r=p}^{N_p}\lambda^p_r H(.,Y^r)\Phi^r = \sum_{r=p}^{N_p}\lambda^p_r \left[H(.,Y^r)-H(.,Y)\right]\Phi^r +  H(.,Y){}^p\!{\Phi}.
\end{align*}
Using the Lipschitz property of $H$ and the uniform $L^2$-bound on $ \nabla \varphi_n(Y^n_{\cdot})$, the first term in the right hand side of the previous equation tends to zero in $\HP{2}$. Then we get
\begin{align*}
  \sum_{r=p}^{N_p}\lambda^p_r H(.,Y^r)\Phi^r  \overset{\HP{2}}{\longrightarrow} H(.,Y)\Phi.
\end{align*}
Finally, we just have to pass to the limit into
$${}^p\!Y_t = \xi+ \int_t^T  \sum_{r=p}^{N_p}\lambda^p_r f(s,Y^r_s,Z^r_s) \ud s -\int_t^T {}^p\!Z_s \ud W_s - \int_t^T \sum_{r=p}^{N_p}\lambda^p_r H(s,Y^r_s)\Phi^r_s \ud s$$
 to conclude the proof of the theorem.

\eproof
%

%\vspace{5pt}
%\begin{Corollary}
%the case $\xi = g(X)$ and the proof is in the next section.
%\end{Corollary}

\subsection{General case}
%\textcolor{magenta}{
%\begin{enumerate}
%\item  donner Resultat principal existence + uniqueness sous forme d'un theoreme
%\item  prouver lemme control\eqref{eq hyp gen Phi}
%\item faire la preuve du theorem en deux étapes 1. unicité, 2. existence via troncation (smooth?) de la condition terminal attention! que les $\xi^n$ appartiennent bien a $\mathfrak{T}_\beta$.
%\end{enumerate}
%}

\begin{Theorem}\label{th main result under SB}
Assume that \textbf{(SB)} holds and $\xi \in \mathfrak{T}_{\mathfrak{B}(\Lambda)}$. There exists a unique solution $(Y,Z,\Phi) \in \SP{2}\times\BP{2}\times\BP{2}$ to \eqref{eq:EDSRORintro}. 
%\textcolor{magenta}{Mettre quelques proprietes de la solution ?}
\end{Theorem}

Before proving our main result, we consider the following lemma which is a key result for the study of Obliquely Reflected BSDEs, as it proves, among other things, the structural condition \eqref{eq hyp gen Phi}. It is the counterpart of Lemma \ref{le control Phi and dist penalised} introduced for the penalised BSDE.

\begin{Lemma}\label{le proof struc cond orbsde}
Assume that \textbf{(SB)} holds. 
Let $(Y,Z,\Phi) \in \SP{2}\times\HP{2}\times\HP{2}$ be a solution to the Obliquely Reflected BSDE \eqref{eq:EDSRORintro}. 
Then, the structural condition \eqref{eq hyp gen Phi} holds true for $(Y,Z,\Phi)$ for some $K > 0$. Moreover, 
there exists $c:=c(\sigma^\xi)$ such that
\begin{align}\label{eq unif control Phi}
\NB{2}{\Phi} \le c\;.
%\esp{\int_0^T |\Phi_s|^2 \ud s} \le c'\esp{|\xi|^2 + \int_0^T|\alpha_s|^2 \ud s}\;.
\end{align}
\end{Lemma}
%\begin{bluetext}
\proof
Applying It\^o's formula to $U_t := \phi(Y_t)$, recall assumption \textbf{(SB)},  we compute that
$\ud U_t = a_t \ud t + b_t \ud W_t $ with
\begin{align*}
%\ud U_t &= a_t \ud t + \beta_t \ud W_t + \ud A_t \quad \text{with}
%\\
a_t&:= \partial \phi(Y_t)\set{-f(t,Y_t,Z_t) +  H(t,Y_t)\Phi_t}+ \frac12 \mathrm{Tr}[\partial^2 \phi(Y_t)Z_tZ^*_t] 
\quad \text{and} \quad b_t := \partial \phi(Y_t)Z_t.
\end{align*}
%and $A$ is an increasing process. 
Using It\^o-Tanaka formula, we obtain
\begin{align*}
%\ud (-U_t )
%&= -\alpha_t \ud t
%- \beta_t \ud W_t \quad \text{ and } 
\quad \ud [-U_t ]^+
= -a_t1_{\set{U_t<0}} \ud t
- b_t 1_{\set{U_t<0}}\ud W_t + \ud L^0_t
\end{align*}
where $L^0$ is the local time at $0$ of the semi-martingale $U$. Taking the difference of the two previous equations, we obtain
\begin{align*}
% -\alpha_t \ud t
%- \beta_t \ud W_t & = -\alpha_t1_{\set{U_t<0}} \ud t
%- \beta_t 1_{\set{U_t<0}}\ud W_t + \ud L^0_t
%\\
0 &= a_t1_{\set{U_t=0}} \ud t
 + b_t 1_{\set{U_t=0}}\ud W_t + \ud L^0_t
\end{align*}
which leads to
$
a_t1_{\set{U_t=0}} \ud t \le 0 
$.
 We then deduce 
\begin{align} \label{eq majo Phi}
 |\Phi_t| \ud t \le \frac1\eta [\partial \phi(Y_t)f(t,Y_t,Z_t)]^+ \ud t,
\end{align}
recall \eqref{eq vector field property} and Remark \ref{rem hyp domain} $i)$. From this, we deduce that \emph{a fortiori} \eqref{eq hyp gen Phi} holds true.
%\\
%Under \textbf{(SB)}, this leads to 
%\textcolor{red}{ On peut ecrire sous (A) puis recuperer uniformite via premier lemme des a priori - faire le lien entre (A) et (SB) ou alors PAS UTILE et on utilise controle fin obtenu dans second a priori estimate sous (SB) à voir...}
\eproof
%\end{bluetext}

\vspace{4pt}
%\begin{bluetext}
\noindent We should notice that in the proof of the above lemma, we obtain a stronger result than the structural condition \eqref{eq hyp gen Phi}. Indeed, we are able to control in \eqref{eq majo Phi} the reflecting process without the conditional expectation appearing in  \eqref{eq hyp gen Phi}. %\textcolor{red}{comment more?}
%\end{bluetext}

\vspace{10pt}
\noindent We now turn to the proof of our main result for this section.

\vspace{4pt}
\noindent \textbf{Proof of Theorem \ref{th main result under SB}}\\
%\begin{bluetext}
1. We first prove uniqueness of the solution. Let $(\eY,\eZ,\ePhi)$ and $(\zY,\zZ,\zPhi)$ be two solutions of \eqref{eq smooth and bounded ORBSDE} in $\SP{2}\times\HP{2}\times\HP{2}$. We first observe that both solutions  satisfies \eqref{eq hyp gen Phi} by application of Lemma \ref{le proof struc cond orbsde} which allows us to invoke Proposition \ref{pr second a priori estimate}. Moreover, both solutions satisfy \eqref{eq prop of iPhi} by definition. Then, a straightforward application of Proposition \ref{pr stability} concludes the proof of this step, noticing that all the terms in the right hand side of \eqref{eq pr stab} are null.
\\
2. We now turn to the existence question.\\
2.a We first approximate $\xi$ by a sequence of bounded random variables $(\xi_N)_{N \ge 1}$. 
Let $(\tau_N)_{N \ge 1}$ be the sequence of stopping time defined by
\begin{align*}
\tau_N := \inf \set{t \ge 0 \,|\, |\cY^\xi_t| \ge N }\wedge T\;,
\end{align*}
and we set $\xi_N := \cY^\xi_{\tau_N}$.
Importantly, we observe that $\xi_N$ satisfies \textbf{(SB)}(i) and it belongs also to the class $\mathfrak{T}_{\mathfrak{B}(\Lambda)}$, indeed 
$\int_0^T |\cZ^{\xi_N}_s|^2 \ud s \le \int_0^T |\cZ^{\xi}_s|^2 \ud s.$ 
For later use, let us also remark that
\begin{align} \label{eq majo unif quant of interest}
\sigma^{\xi^N} \le \sigma^\xi\;,\; \text{ for all } N\ge 1\;,
\end{align}
recall \eqref{eq de sup Yxi expo int}.
Moreover, since
\begin{align*}
\xi_N \rightarrow \xi \;\P-a.s. \quad \text{ and }\quad |\xi_N-\xi|\le 2 \sup_{t \in [0,T]}|\cY^\xi_t|\,,
\end{align*}
we have that by the dominated convergence theorem, recall Remark \ref{re comment hyp xi} (i), $\xi_N \rightarrow \xi$ in $\LP{q}$, for any $q \ge 1$. 
\\
2.b Applying Proposition \ref{pr smooth and bounded exist orbsde}, we introduce a sequence of Obliquely RBSDEs, $(Y^N,Z^N,\Phi^N)$ with terminal condition $\xi_N$.
We now show that $(Y^N,Z^N)$ is a Cauchy sequence in $\SP{2}\times\HP{2}$. First, we apply the stability estimate given in Proposition \ref{pr stability}: for $N, P \ge 1$, we have
\begin{align*}
\sup_{t \in [0,T]} \esp{|Y^N_t-Y^{P}_t|^2} + \NH{2}{Z^N-Z^{P}}^2
&\le
C\esp{\Gamma_T^{N,P}|\xi^N-\xi^{P}|^2}\;,
\end{align*}
with $\Gamma^{N,P}$  such that for some $p>1$ and $C>0$,
\begin{align*}
\esp{|\Gamma^{N,P}_T|^p} \le C\;,
\end{align*}
where importantly $p$ and $C$ do not depend on $(N,P)$, recall \eqref{eq majo unif quant of interest}. Using H\"older inequality, we then obtain
\begin{align*}
\sup_{t \in [0,T]} \esp{|Y^N_t-Y^{P}_t|^2} + \NH{2}{Z^N-Z^{P}}^2
&\le C\NL{2q}{\xi^N-\xi^{P}}^2\;.
\end{align*} 
Following classical arguments, see Step 1.a in the proof of Proposition \ref{pr smooth and bounded exist orbsde}, we compute also
\begin{align*}
\NS{2}{Y^N-Y^P}^2 \le C\esp{\int_0^T|Y^N_t-Y^P_t|\left(|\Phi^N_t|+|\Phi^P_t| \right)\ud t}
+C\NH{2}{Z^N-Z^P}^2 \;.
\end{align*}
Applying Cauchy-Schwarz inequality, and combining Lemma \ref{le proof struc cond orbsde} and \eqref{eq majo unif quant of interest}, we  get
\begin{align*}
\NS{2}{Y^N-Y^P}^2  \le C \left( \NH{2}{Y^N-Y^P} +  \NH{2}{Z^N-Z^P}^2 \right)\;.
\end{align*}
Eventually, we obtain
\begin{align*}
\NS{2}{Y^N-Y^P}^2  +  \NH{2}{Z^N-Z^P}^2
\le C\left(\NL{2q}{\xi^N-\xi^{P}} + \NL{2q}{\xi^N-\xi^{P}}^2 \right)\;.
\end{align*}
From the conclusion of Step 1. we deduce the Cauchy property of the sequence $(Y^N,Z^N)$ and we denote $(Y,Z)$ its limit. The proof is then concluded following the same arguments as in step 2 of Proposition \ref{pr smooth and bounded exist orbsde}, once  observed that by Lemma
\ref{le proof struc cond orbsde},
$$ \esp{\int_0^T |\Phi_s^N|^2 \ud s } \le C\;,$$
where again $C$ does not depend on $N$ from \eqref{eq majo unif quant of interest}.
\eproof
%\end{bluetext}

%\paragraph{Existence and uniqueness result...}

%\input{sec3-secondpart}

%\input{ultimate}

%\input{bornitude}
%!TEX root = main.tex
\section{A general existence result in the Markovian framework}
\label{se markov}

In this section, we introduce a Markovian framework: for all $(t,x) \in [0,T] \times \R^q$, we denote $(X_s^{t,x})_{s \in [0,T]}$ the solution of the SDE
\begin{align}
\label{SDE}
 \ud X_s &= b(s,X_s)\ud s + \sigma(s,X_s)  \ud W_s, \quad s \in [t,T],\\
 \nonumber X_s &= x, \quad s \in [0,t],
\end{align}
and we consider the following Markovian reflected BSDE:
\begin{equation}
 \label{BSDE reflected}
 \left\{ \begin{aligned} & Y_t = g(X_T^{0,a}) + \int_t^T f(s,X_s^{0,a},Y_s,Z_s) \ud s -\int_t^T Z_s \ud W_s - \int_t^T H(s,X_s^{0,a},Y_s,Z_s)\Phi_s ds, \\
         & Y_t \in \bar{\cD}, \quad \Phi_t \in\partial \varphi(Y_t),\quad 0 \leqslant t \leqslant T,\quad 
          \int_0^T  \1_{\set{Y_t \notin \partial \cD}}  |\Phi_t| \ud t= 0.
        \end{aligned}
 \right.
\end{equation}
The main goal of this section is to prove an existence result for the above reflected BSDE when $H$ is only continuous, compare with assumption \textbf{(SB)}. We also discuss the case of discontinuous $H$ and the difficulty arising for uniqueness in this setting.

\subsection{Continuous oblique direction of reflection}
\label{subse markov continuous}
We now introduce the main setting for this part. The set of assumption below echoes assumption \textbf{(A)} introduced in Section \ref{subse framework} but in a Markovian setting.

\paragraph{Assumption (AM)}
$ $ %\textcolor{red}{dire que H et f verifient A}
\begin{enumerate}[i)]
 %\item The boundary $\partial \cD$ of the domain can be written as the following partition 
 %$$\partial \cD =  \mathcal{O} \cup \bigcup_{i \in I} \mathcal{C}_i$$ 
 %where $I$ is a finite set, $\mathcal{O}$ is an open set of $\cD$ and for all $i \in I$, $\mathcal{C}_i$ is a closed set of $\cD$.
 \item $b: [0,T] \times \R^q \rightarrow \R^q$ and $\sigma: [0,T] \times \R^q \rightarrow \R^{q \times k}$ are measurable functions satisfying linear growth condition and uniform Lipschitz condition with respect to the space variable namely
 \begin{align*}
 |b(t,x)| + |\sigma(t,x)| \le L(1+|x|) \text{ and }  |b(t,x)-b(t,y)| + |\sigma(t,x)-\sigma(t,y)| \le L|x-y|\;,
 \end{align*}
 for some $L>0$ and all $(t,x,y)\in [0,T]\times\R^q\times\R^q$.
 \item $g: \R^q \rightarrow \bar{\cD}$ is a measurable function and there exists $p\in \R^+$ such that for any $x \in \R^q$,
 $$|g(x)| \leq L(1+|x|^p).$$
 \item $f: [0,T] \times \R^q \times \R^d \times \R^{d\times k}\rightarrow \R^{d}$ is a measurable function satisfying: there exists $p \in \R^+$ such that, for any $(t,x,y,z) \in [0,T] \times \R^q \times \R^d \times \R^{d\times k}$, we have
 $$|f(t,x,y,z)| \leq L(1+|x|^p+|y|+|z|),$$
 and, for all $(t,x) \in [0,T] \times \R^d$, $f(t,x,.,.)$ is continuous on $\R^d \times \R^{d\times k}$.
 \item $H: [0,T] \times \R^q \times  \R^d \times \R^{d\times k}\rightarrow \R^{d \times d}$ is a measurable function. There exists $\eta > 0$  such that, for all $(t,x,y,z) \in [0,T] \times \R^q \times \mathbb{R}^d  \times \R^{d \times k}$%, %$H(t,x,y,z)$ is a symetric matrix 
 %\textcolor{red}{(OK:on a pas besoin de symétrie a priori...)}, 
 \begin{align}
  \label{inegalite H a}
   H(t,x,y,z) \upsilon \cdot \upsilon \geq \eta ,\quad \forall \upsilon \in \mathfrak{n}(\mathfrak{P}(y)) ,% |y-\mathfrak{P}(y)|^2,}%\quad \forall u \in \R^q, 
 \end{align}
 %\textcolor{red}{on en a pas besoin dans toutes les directions et on a pas besoin de l'inversibilite globale...}
 $$\text{and} \quad |H(t,x,y,z)| \leq L.$$
%and
% $$H(t,x,y,z)=0 \quad \text{when } y\in \cD.$$
%and  there exists an open set $\mathcal{V}$ of $\mathbb{R}^d$ such that $\mathcal{O} \subset \mathcal{V}$ and for all $(t,x) \in [0,T] \times \R^d$, $H(t,x,.,.)$ is uniformly continuous on $\mathcal{V} \times \R^{d\times k}$.
\item Let $\mathcal{X} = \{\mu(t,x;s,dy),x \in \R^q \textrm{ and } 0\leq t \leq s \leq T\}$ be the family of laws of $X^{t,x}$ on $\R^q$, i.e., the measures such that $\forall A \in \mathcal{B}(\R^q)$, $\mu(t,x;s,A) = \mathbb{P}(X_s^{t,x} \in A)$. 
%We will assume that $\mathcal{X}$ satisfies a $L^2$-domination condition.
%\paragraph{Assumption (HDom)}
For any $t \in [0,T)$, for any $\mu(0,a;t,dy)$-almost every $x \in \R^q$, and any $\delta \in ]0,T-t]$, there exists an application $\phi_{t,x}: [t,T]\times \R^d \rightarrow \R^+$ such that:
\begin{enumerate}
 \item $\forall k \geq 1$, $\phi_{t,x} \in L^2([t+\delta,T] \times [-k,k]^q; \mu(0,a;s,dy)ds)$,
 \item $\mu(t,x;s,dy)ds = \phi_{t,x}(s,y)\mu(0,a;s,dy)ds$ on $[t+\delta,T] \times \R^q$.
\end{enumerate}
\item  For all $(t,x) \in [0,T] \times \R^d$, $H(t,x,.,.)$ is continuous on $\R^d \times \overline{\cD}$.
\end{enumerate}

\begin{Remark}
\begin{enumerate}[i)]
 \item We observe that $H(t,X,\cdot)$ and $f(t,X,\cdot)$ satisfy assumption \textbf{(A)}. Thus we will use in the sequel the \emph{a priori} estimates obtained in Section \ref{subse a priori estimates}.
 %\textcolor{red}{dire que H et f verifient A}
 \item Remark \ref{re extension H nonverysmooth} applies for $H$ which is continuous in this context.
% In applications, for all $(t,x,z) \in [0,T] \times \mathbb{R}^k \times \mathbb{R}^{d\times k}$ $H(t,x,.,z)$ is usually specified only on the boundary $\partial \cD$. Nevertheless, if $H(t,x,.,z)$ is a continuous and bounded fonction on $\partial \cD$ it is possible to extend it to a continuous and bounded function on $\cD$. Indeed, $\cD$ is homeomorph to a set $S$ which is a half plane of $\mathbb{R}^d$ or $\mathbb{R}^r \times B^{d -r}$, with $0 \leqslant r \leqslant d$. Moreover, the boundary of $\cD$ is sent to the boundary of $S$. Then we just have to remark that the extension of $H(t,x,.,z)$ is easy when $\cD=S$. 
%
 \item The $\mathscr{L}^2$-domination condition {\bf(AM)}(v) was initially introduced in \cite{Hamadene-Lepeltier-Peng-97}. We refer to \cite{Hamadene-Lepeltier-Peng-97,DeAngelis-Ferrari-Hamadene-17} for examples of assumptions on coefficients of the SDE \eqref{SDE} under which {\bf(AM)}(v) is true.%\textcolor{magenta}{
%$L^2$-domination condition as  in \cite{Hamadene-Lepeltier-Peng-97}.
%}
\end{enumerate}
\end{Remark}

%%%%%%%%%%%%%%%%%%%%%%%%

\begin{Theorem}
\label{theorem existence}
 Assume \textbf{(AM)}. Then, there exists a solution $(Y,Z) \in \mathscr{S}^2 \times \mathscr{H}^2$ to \eqref{BSDE reflected}. Moreover,
 the following Markovian representation holds true: \\There exist $u : [0,T] \times \R^q \rightarrow \R^d$ and $v : [0,T] \times \R^q \rightarrow \R^{d \times k}$  measurable functions such that
 \begin{equation*}
  Y_t=u(t,X_t^{0,a}) \quad \textrm{and} \quad Z_t=v(t,X_t^{0,a}),
 \end{equation*}
and, for some $c>0$, for all $(t,x) \in [0,T] \times \R^q$,
$$|u(t,x)| \leq c(1+|x|^p).$$
\end{Theorem}

By choosing properly the function $H$ we can obtain the following corollary.

\begin{Corollary}
 \label{corollaire existence EDSR switching}
 Let us consider the following obliquely reflected Markovian BSDE
 \begin{equation}\label{eqBSDE swithching general}\qquad
\begin{cases}
\displaystyle
Y_t = g(X_T^{0,a})+\int_t^T f(s,X_s^{0,a},Y _s, Z _s)\,\ud s-\int
_t^T Z_s \,\ud W_s  + \int_t^T \Psi_s\ud s, &\quad 0 \le t \le
T,\vspace*{2pt}\cr
\displaystyle Y_t^\ell \ge \max_{j \in \cI} \set{Y_t^j -c^{\ell j}} ,  &\quad
0\le t\le T,\; \ell \in \cI, \vspace*{2pt}\cr
\displaystyle\int_0^T\Bigl[Y _t^\ell-\max_{\textcolor{black}{j\in \cI \setminus \set{\ell} } }\{Y ^j_t -
c^{\ell j}\}\Bigr]
\,\Psi_t^\ell \ud t=0,
&\quad \ell \in
\cI,
\end{cases}
\end{equation}
where $\cI:=\{1,\ldots,d\}$ and the switching costs $(c^{ij})_{i,j\in\cI}$ satisfy the following structure condition
 \begin{equation}
  \label{condition de structure c}
  \begin{cases}
   c^{ii}=0, & \textrm{for } 1 \le i \le d;\\
   % c^{ij} >0, & \textrm{for } 1 \le i,j \le d \textrm{ with } i \neq j;\\
   \{c^{ij}+c^{jl} - c^{il}\} >0, & \textrm{for } 1 \le i,j \le d \textrm{ with } i \neq j, j \neq l.
  \end{cases}
 \end{equation}
 We assume that assumption \textbf{(AM)} is in force. Then there exists a solution $(Y,Z,\Psi) \in \SP{2} \times \HP{2} \times \HP{2}$ to \eqref{eqBSDE swithching general} . Moreover
 the following Markovian representation holds true: There exist two measurable functions $u : [0,T] \times \R^q \rightarrow \R^d$ and $v : [0,T] \times \R^q \rightarrow \R^{d \times k}$ such that
 \begin{equation*}
  Y_t=u(t,X_t^{0,x}) \quad \textrm{and} \quad Z_t=v(t,X_t^{0,x}),
 \end{equation*}
and, for some $c>0$, for all $(t,x) \in [0,T] \times \R^q$,
$$|u(t,x)| \leq c(1+|x|^p).$$
\end{Corollary}
\begin{Remark}
 %\textcolor{red}{difference avec resultats precedents biblio... }
 The main novelty here is the dependence  of the generator on the whole $z$ (as in the concomitant article \cite{DeAngelis-Ferrari-Hamadene-17}) which extend the result of \cite{Hu-Tang-10,Hamadene-Zhang-10,Chassagneux-Elie-Kharroubi-11} and the possibility to consider negative switching costs. We refer to \cite{Martyr-16} and references inside for a recent work dealing with switching problems with signed switching costs. Our result only cover the case of constant switching costs due to a priori estimates obtained previously in the framework of a deterministic domain $\cD$. Nevertheless our approach might be adapted to treat random domains and then tackle the problem of switched BSDEs with random signed switching costs.
\end{Remark}

Before giving the proof of Theorem \ref{theorem existence} and Corollary \ref{corollaire existence EDSR switching}, we start by considering an approximation of \eqref{BSDE reflected}. 
%We define, for all $(t,x,y,z) \in [0,T] \times \mathbb{R}^k \times \mathbb{R}^d \times \mathbb{R}^{d\times k}$,
%$$\bar{H}(t,x,y,z) = H(t,x,\mathfrak{P}(y),z).\text{ \textcolor{red}{deja def ? simplement etendre H des le depart ?}}$$
%\begin{cases}
%                      H(t,x,\mathfrak{P}(y),z) &\text{if } y \notin \cD,\\
%                      0 & \text{if } y \in \cD.
%                     \end{cases}$$
Let $\theta$ be an element of $C^\infty (\mathbb{R}^{d+d\times k},\mathbb{R}^+)$ with compact support and satisfying 
$$\int_{\mathbb{R}^{d+d\times k}} \theta(y,z)\ud y\ud z=1.$$
For all $n \in \mathbb{N}$ and $(t,x,y,z) \in [0,T] \times \mathbb{R}^q \times \mathbb{R}^d \times \mathbb{R}^{d\times k}$ we set
\begin{align*}
 f_n(t,x,y,z) &=  \int_{\mathbb{R}^{d+d\times k}} n^2 f(t,x,y,z) \theta(n(y-u),n(z-v)) \ud u \ud v\\
 H_n(t,x,y,z) &= \int_{\mathbb{R}^{d+d\times k}} n^2 H(t,x,y,z) \theta(n(y-u),n(z-v)) \ud u \ud v.
\end{align*}
By classical convolution arguments %\textcolor{red}{(donner explicitement convolution. Pour H necessaire prolonger par $0$ sur $\cD$ puis de convoluer avec une fonction a support compact)} 
functions $(f_n)_{n \in \mathbb{N}}$ and $(H_n)_{n \in \mathbb{N}}$ satisfy following properties.
\begin{Lemma}\label{le smooth f and H} Assume \textbf{(AM)}.
\begin{enumerate}[i)]
 \item $f_n : [0,T] \times \mathbb{R}^{q} \times \mathbb{R}^d \times \mathbb{R}^{d \times k} \rightarrow \mathbb{R}^d$ and $H_n : [0,T] \times \mathbb{R}^{q} \times \mathbb{R}^d \times \mathbb{R}^{d \times k} \rightarrow \mathbb{R}^{d\times d}$ are measurable and uniformly Lipschitz functions with respect to $(y,z)$.
 \item $|f_n(t,x,y,z)| \le L(1+|x|^p+ |y|+|z|)$ and $|H_n(t,x,y,z)| \le L $ for all $(t,x,y,z) \in [0,T] \times \mathbb{R}^{q} \times \mathbb{R}^{d}  \times \mathbb{R}^{d \times k}$.
 \item For all $(t,x) \in [0,T] \times \mathbb{R}^q$ and $\cK$ a compact subset of $\mathbb{R}^d\times \mathbb{R}^{d \times k}$
 $$\sup_{(y,z) \in \cK} |f_n(t,x,y,z)-f(t,x,y,z)| + \sup_{(y,z) \in \cK} |H_n(t,x,y,z)-H(t,x,y,z)| \stackrel{n \rightarrow +\infty}{\longrightarrow} 0.$$
 %\item $H_n : [0,T] \times \mathbb{R}^{q} \times \bar{\cD} \times \mathbb{R}^{d \times k} \rightarrow \mathbb{R}^{d\times d}$ is measurable and uniformly Lipschitz with respect to $(y,z)$.
 %\item $|H_n(t,x,y,z)| \le L$ for all $(t,x,y,z) \in [0,T] \times \mathbb{R}^{q} \times \bar{\cD} \times \mathbb{R}^{d \times k}$.
 %\item $H_n$ converges to $H$ a.e.
 %\item For any $(t,x) \in [0,T] \times \mathbb{R}^q$
 %$$\sup_{(y,z) \in \mathcal{V} \times \mathbb{R}^{d\times k}} |H_n(t,x,y,z)-H(t,x,y,z)| \stackrel{n \rightarrow +\infty}{\rightarrow} 0.$$
\end{enumerate}
%\textcolor{red}{on a pas besoin de la convergence de $H_n$ vers $H$, car on veut juste caractériser la direction de projection}.
\end{Lemma}

\vspace{4pt}
\noindent For any $n \in \mathbb{N}$, we then consider the following BSDE
\begin{align}
 \nonumber
 Y^n_t =& g(X_T^{0,a}) + \int_t^T f_n(s,X_s^{0,a},Y^n_s,Z^n_s) \ud s\\
 &-\int_t^T Z^n_s \ud W_s - \int_t^T H_n(s,X_s^{0,a},Y^n_s,Z^n_s)\nabla \varphi_n (Y^n_s)\ud s , \quad t \in [0,T] \label{eq de penalized}
\end{align}
where $\varphi_n$ is defined in \eqref{de varphi n} {with $M$ fixed to an arbitrary value}. Note that, in this section, for the reader's convenience, we write simply $\varphi_n$ instead  of $\varphi_n^M$.

\begin{Lemma}
 There exists a unique solution to \eqref{eq de penalized} in $\mathscr{S}^2\times\mathscr{H}^2$. Moreover, there is a Markovian representation for this solution: for all $n \in \mathbb{N}$, there exist $u_n:[0,T]\times \R^q \rightarrow \R^d$ and $v_n: [0,T] \times \R^q \rightarrow \R^{d\times k}$  measurable functions satisfying
\begin{align}
 Y^{n}_t = u_n(t,X^{0,a}_T) \text{ and } Z^n_t = v_n(t,X^{0,a}_T).
\end{align}
Moreover, for all $(t,x) \in [0,T] \times \R^q$, $(u_n(s,X_s^{t,x}),v_n(s,X_s^{t,x}))_{s \in [t,T]}$ is the unique solution in $\mathscr{S}^2 \times \mathscr{H}^2$ of the BSDE 
\begin{align}
 \nonumber Y^{n,t,x}_s =& g(X_T^{t,x}) + \int_s^T f_n(r,X_r^{t,x},Y^{n,t,x}_r,Z^{n,t,x}_r) \ud r -\int_s^T Z^{n,t,x}_r \ud W_r\\
 \label{eq de penalized t x}
 &- \int_s^T H_n(r,X_r^{t,x},Y^{n,t,x}_r,Z^{n,t,x}_r)\nabla \varphi_n (Y^{n,t,x}_r)\ud r \quad s \in [t,T].
\end{align}
\end{Lemma}
\proof
We use the same arguments as in the proof of Lemma \ref{le:euEDSRpenalise}:
Since $H_n$ and $\nabla \varphi_n$ are two Lipschitz bounded functions (with respect to $y$ and $z$),  the penalised BSDE \eqref{eq de penalized} has a Lipschitz driver and the classical theory then applies to get the existence, uniqueness and representation result.
\eproof
 
\vspace{4pt} 
By applying Lemma \ref{le first basic a priori} and Lemma \ref{le control Phi and dist penalised}, we obtain the following estimates for $(Y^{n,t,x},Z^{n,t,x})$.

\begin{Proposition}\label{pr apriori}
 For all $(t,x) \in [0,T] \times \R^q$, we have
 \begin{align*}
\sup_{t \le s \le T} \mathbb{E}\left[|Y^{n,t,x}_s|^2+ \varphi_n(Y^{n,t,x}_s)\right] 
+\mathbb{E}\left[\int_t^T|Z^{n,t,x}_s|^2 \ud s+\int_t^T|\nabla \varphi_n(Y^{n,t,x}_s)|^2\ud s\right] \le C(1+|x|^{2p})\;.
 \end{align*}
\end{Proposition}
\noindent In particular, Proposition \ref{pr apriori} yields that, for some $c>0$,%for all $n \in \mathbb{N}$ and $(t,x) \in [0,T] \times \R^q$,
$$|u_n(t,x)| \leq c(1+|x|^p), \quad \forall n \in \mathbb{N}, \quad \forall(t,x) \in [0,T] \times \R^q.$$% \text{for some } c>0.$$

\noindent We now turn to the proof of the main result for this section.

\vspace{5pt}
\noindent \textbf{Proof of Theorem \ref{theorem existence}}\\
The proof follows mainly from arguments in \cite{Hamadene-Lepeltier-Peng-97}. Some extra work is required to identify the reflecting process properly.

\noindent 1. Define,
\begin{align*}
 F_n(t,x) = f_n(t,x,u_n(t,x),v_n(t,x)),\quad G_n(t,x) = H_n(t,x,u_n(t,x),v_n(t,x))\nabla \varphi_n(u_n(t,x))\;,
\end{align*}
and 
$$\mathfrak{F}_n := F_n - G_n,$$
 we compute
\begin{align*}
\int_{\R^q}\int_0^T |\mathfrak{F}_n(s,y)|^2 \mu(0,a;s,\ud y) \ud s &= \esp{\int_0^T |\mathfrak{F}_n(s,X_s^{0,a})|^2 \ud s}
\\
&\le \esp{\int_0^T C(1 + |X_s^{0,a}|^{2p} + |Y^n_s|^2 + |Z^n_s|^2
+ |\nabla \varphi_n(Y^n_s)|^2
) \ud s}
\\
& \le C,
\end{align*}
by using Proposition \ref{pr apriori}. Thus we get
 %\begin{align*}
  $\mathfrak{F}_n \rightharpoonup \mathfrak{F}$
 %\end{align*}
 in $L^2([0,T]\times \R^q;\mu(0,a;s,\ud x) \ud s)$, up to a subsequence.

 \vspace{2pt}
 2. We now show that $(u_n(t,x))_{n \in \mathbb{N}}$ is a Cauchy sequence in $\R^d$ for all $t\in [0,T]$ and for $\mu(0,a;t,\ud x)$-almost every $x \in \R^q$. 
 %First of all, when $t=0$, $(u_n(0,a))_{n \in \mathbb{N}}$ is a bounded sequence, so it converges up to a subsequence. For simplicity, we use the same notation for this subsequence. Now we assume that $t>0$.
 When $t=T$ the sequence is constant and the result is obvious. When $t<T$, $x\in \R^q$ and $\delta \in (0,T-t]$,
 %$(t,x) \in (0,T]\times \R^q$. 
 we compute
 \begin{align*}
  |u_n(t,x)-u_m(t,x)|& = \left|\esp{\int_t^T\left(\mathfrak{F}_n(s,X^{t,x}_s) - 
		               \mathfrak{F}_m(s,X^{t,x}_s)\right)\ud s  }\right|
		   \\
&\le \esp{\int_t^{t+\delta}|\mathfrak{F}_n(s,X^{t,x}_s) - 
		               \mathfrak{F}_m(s,X^{t,x}_s)| \ud s} =:A_1
\\ 
& \;\; + \esp{\int_{t+\delta}^T
|\mathfrak{F}_n(s,X^{t,x}_s) - \mathfrak{F}_m(s,X^{t,x}_s)| \1_{\set{|X_s^{t,x}| \ge \kappa}}\ud s
} =:A_2 
\\
&\;\; +
 \left|\esp{\int_{t+\delta}^T
\left(\mathfrak{F}_n(s,X^{t,x}_s) - \mathfrak{F}_m(s,X^{t,x}_s)\right)  \1_{\set{|X_s^{t,x}| < \kappa}}\ud s
}\right| =:A_3 \;.
 \end{align*}
For the first two terms, we easily get
\begin{align*}
 A_1 &\le \delta^\frac12\esp{\int_t^{t+\delta}|\mathfrak{F}_n(s,X^{t,x}_s) - 
		               \mathfrak{F}_m(s,X^{t,x}_s)|^2 \ud s}^\frac12 \le C(1+|x|^p)\delta^\frac12\;,
	\\
A_2 &\le C \kappa^{-\frac12}\esp{\int_{t+\delta}^T |X_s^{t,x}|\ud s}\esp{\int_{t+\delta}^T|\mathfrak{F}_n(s,X^{t,x}_s) - 
		               \mathfrak{F}_m(s,X^{t,x}_s)|^2 \ud s}^\frac12 \le C(1+|x|^{p+1}) \kappa^{-\frac12} \;,
\end{align*}
where $C$ is a constant that does not depend on $n$ nor $m$.
For the third term, we have
\begin{align*}
 A_3 &= \left| \int_{\R^q}\int_{t+\delta}^T \left( 
 \mathfrak{F}_n(s,y) -   \mathfrak{F}_m(s,y)
 \right) \1_{\set{|y|\le \kappa}} \mu(t,x;s,\ud y)\ud s \right|
 \\
 &=\left|\int_{\R^q}\int_{t+\delta}^T \left( 
 \mathfrak{F}_n(s,y) -  \mathfrak{F}_m(s,y)
 \right) \1_{\set{|y|\le \kappa}} \phi_{t,x}(s,y)\mu(0,a;s,\ud y)\ud s \right|
\end{align*}
for $\mu(0,a;s,\ud x)$-almost every $x \in \R^q$,
%for $t>0$, $x \in \R^q$ \textcolor{red}{il faut mettre des $\mu(...)$ pp},
where we used the $L^2$-domination assumption.
By weak convergence,
$A_3 \rightarrow 0$ when $n,m \rightarrow \infty$.
Thus, by taking $\delta \rightarrow 0$ and $\kappa \rightarrow +\infty$ we show that for all $t\in [0,T]$ and for $\mu(0,a;t,\ud x)$-almost every $x \in \R^q$, $(u_n(t,x))_{n \in \mathbb{N}}$ is a Cauchy sequence. %For $t=0$, $u_n(0,a)$ is bounded, so it converges up to a subsequence. Eventually, we can define
So, there exists a Borelian application $u : [0,T] \times \R^q \rightarrow \R^d$ such that for all $t\in [0,T]$, for $\mu(0,a;t,\ud x)$-almost every $x \in \R^q$,
\begin{align}
u(t,x) = \lim_{n\infty}u_n(t,x).
\end{align}
\\
We straightforwardly get, for all $t \in [0,T]$,
\begin{align*}
 Y^n_t = u_n(t, X^{0,a}_t) \rightarrow u(t,X^{0,a}_t):=Y_t, \quad \textrm{a.s.}
\end{align*}
and , observing that $|Y^n_t| \le C(1+|X_t^{0,a}|^p)$, we obtain via the dominated convergence theorem, $Y^n_t \rightarrow Y_t$ in $L^2([0,T]\times \Omega,dt \otimes d\P)$.

\vspace{2pt}
3. We can easily prove that the process $Y$ lives in the convex set $\bar{\cD}$. Indeed, we have, recalling \eqref{eq expression  varphi},
\begin{align*}
 \sup_{0 \le s \le T} \mathbb{E}\left[ \varphi_1 (Y_s)\right] &\le  \sup_{0 \le s \le T} \mathbb{E}\left[ |\varphi_1 (Y_s) -\varphi_1 (Y_s^n)|\right]+\frac{1}{n}\sup_{0 \le s \le T} \mathbb{E}\left[ \varphi_n (Y_s^n)\right]\\
 & \le M \sup_{0 \le s \le T} \mathbb{E}\left[ |Y_s - Y_s^n|\right]+\frac{C}{n} \stackrel{n\rightarrow + \infty}{\longrightarrow} 0,
\end{align*}
where we used Proposition \ref{pr apriori}, the fact that $\varphi_1$ is a $M$-Lipschitz function and the convergence of $(Y^n)_{n \in \mathbb{N}}$. Then, for all $s \in [0,T]$, $d(Y_s,\cD)=0$ a.s. and so $Y_s \in \bar{\cD}$ a.s.

\vspace{2pt}
\noindent 4. We now show that $(Z^n)_n=((v_n(t,X_t^{0,a}))_{t \in [0,T]})_n$ is a Cauchy sequence in $L^2([0,T]\times \Omega,dt \otimes d\P)$. For $n,m \ge 1$, we compute, applying It\^o's formula,
\begin{align*}
 \esp{\int_0^T|Z^n_s -Z_s^m|^2 \ud s} &\le 2 \esp{\int_0^T (Y_t^n - Y_t^m)\left(\mathfrak{F}_n(t,X^{0,a}_t) - \mathfrak{F}_m(t,X^{0,a}_t)\right)\ud t} 
 \\
 &\le C  \esp{\int_0^T|Y^n_t-Y^m_t|^2 \ud t} ^\frac12\;,
\end{align*}
which goes to $0$ as $n,m \rightarrow \infty$. 
We denote by $Z$ the limit. From now on, we work with the progressively measurable version of $(Y,Z)$.

\vspace{2pt}
\noindent 5.a In the last step we have to prove that $(Y,Z)$ is a solution to BSDE \eqref{BSDE reflected}. We start by studying the convergence of the generator. Firstly, we compute, for all $\kappa>0$,
\begin{align*}
 &\esp{\int_0^T |f_n(s,X_s^{0,a},Y_s^n,Z_s^n)-f(s,X_s^{0,a},Y_s,Z_s)|\ud s}\\
 \le & \esp{\int_0^T |f_n(s,X_s^{0,a},Y_s^n,Z_s^n)-f(s,X_s^{0,a},Y_s^n,Z_s^n)|\1_{\set{|Y^n_s|+|Z_s^n| \le \kappa}}\ud s}=:B_1\\
 & + \esp{\int_0^T |f_n(s,X_s^{0,a},Y_s^n,Z_s^n)-f(s,X_s^{0,a},Y_s^n,Z_s^n)|\1_{\set{|Y^n_s|+|Z_s^n| > \kappa}}\ud s}=:B_2\\
 & +\esp{\int_0^T |f(s,X_s^{0,a},Y_s^n,Z_s^n)-f(s,X_s^{0,a},Y_s,Z_s)|\ud s}=:B_3.
\end{align*}
Since $f$ and $f_n$ have a linear growth that does not depend on $n$, and $(Y^n,Z^n)$ is uniformly bounded in $L^2([0,T] \times \Omega, dt \otimes dP)$, we get, by using Markov inequality,
\begin{align*}
 B_2 \le & \frac{C}{\kappa}.
\end{align*}
Moreover, we also get
\begin{align*}
 B_3 \leq \esp{ \int_0^T \frac{|f(s,X_s^{0,a},Y_s^n,Z_s^n)-f(s,X_s^{0,a},Y_s,Z_s)|^2}{(1+|Y_s^n|+|Z_s^n|)^2} \ud s}^{1/2}\esp{\int_0^T (1+|Y_s^n|+|Z_s^n|)^2 \ud s}^{1/2}
\end{align*}
and thus, the dominated convergence theorem yields that $B_3$ converges to $0$ as $n\rightarrow +\infty$.

\noindent We now study the first term $B_1$. We have, for all $\kappa>0$, $s \in [0,T]$,
$$|f_n(s,X_s^{0,a},Y^n_s,Z^n_s)-f(s,X_s^{0,a},Y^n_s,Z^n_s)|\1_{\set{|Y^n_s|+|Z_s^n| \le \kappa}} \le C(1+2\kappa+|X_s^{0,a}|),$$
and
\begin{align*}
 &|f_n(s,X_s^{0,a},Y^n_s,Z^n_s)-f(s,X_s^{0,a},Y^n_s,Z^n_s)|\1_{\set{|Y^n_s|+|Z_s^n| \le \kappa}}\\
 \le & \sup_{(y,z), \, |y|+|z| \le \kappa} |f_n(s,X_s^{0,a},y,z)-f(s,X_s^{0,a},y,z)|.
\end{align*}
Thanks to Lemma \ref{le smooth f and H}(iii) we can assert that the second term of the last inequality converges to $0$ and then, by applying the dominated convergence theorem, $B_1$ converges also to $0$. By taking $\kappa \rightarrow +\infty$, it follows that $(f_n(t, X_t^{0,a}, Y^n_t,Z^n_t))_{t \in [0,T]}$ converges to $(f(t, X_t^{0,a}, Y_t,Z_t))_{t \in [0,T]}$ in $L^1([0,T] \times \Omega, dt \otimes dP)$.
\\
%\begin{bluetext}

\noindent 5.b Finally we study the reflecting term.
Since
$$  \esp{\int_0^T|\nabla \varphi_n(Y^n_s)|^2 \ud s} \le C,$$
we have, up to a subsequence, the following weak $L^2([0,T]\times\Omega)$-convergence:
%$$\nabla \varphi_n(Y^n_\cdot ) \rightharpoonup \Psi, \quad \text{when } n \rightarrow + \infty,$$
%and
$$\nabla \varphi_n(Y^n_\cdot ) \rightharpoonup \Phi, \quad \text{when } n \rightarrow + \infty,$$
and we can follow step 2.a in the proof of Proposition \ref{pr smooth and bounded exist orbsde} to obtain
\begin{align*}
\Phi \in \partial \varphi(Y) \; \text{ and } \;  \int_0^T  \1_{\set{Y_t \notin \partial \cD}}  |\Phi_t|\ud t = 0 \;,
\end{align*}
which fully characterize $\Phi$. We now follow step 2.b in the proof of Proposition \ref{pr smooth and bounded exist orbsde}.
Using Mazur's Lemma, we know that there exists a convex combination of $(\Phi^n)_{n \in \mathbb{N}}:=(\nabla \varphi_n(Y^n))_{n \in \mathbb{N}}$ converging strongly in $L^2([0,T]\times\Omega)$, namely
\begin{align*}
 {}^p\!{\Phi} := \sum_{r=p}^{N_p}\lambda^p_r \Phi^r \stackrel{p \rightarrow \infty}{\rightarrow} \Phi,
\end{align*}
where $\lambda^p_r \ge 0$ for all $p \in \mathbb{N}$ and $p \le r \le N_p$, and $\sum_{r=p}^{N_p}\lambda^p_r =1$.
Let us observe that by strong convergence, the following combination
\begin{align*}
 ({}^p\!Y,{}^p\!Z) := \sum_{r=p}^{N_p}\lambda^p_r (Y^r,Z^r)
\end{align*}
still converges to $(Y,Z)$ in $\SP{2} \times \HP{2}$ and, by strong convergence again, 
\begin{align*}
 \sum_{r=p}^{N_p}\lambda^p_r f_r(\cdot,X^{0,a},Y^r,Z^r) &\overset{L^1([0,T] \times \Omega, dt \otimes dP)}{\longrightarrow}
  f(\cdot,X^{0,a},Y,Z)\; \text{ and } \; \int_0^t {}^p\!Z_s \ud W_s \overset{\LP{2}}{\longrightarrow} \int_0^t Z_s \ud W_s \;.
\end{align*}
Moreover, we remark that, for all $t \le T$,
\begin{eqnarray}
\cE^p &:= &\int_0^t  \sum_{r=p}^{N_p}\lambda^p_r  H_r(s,X^{0,a}_s,Y^r_s,Z^r_s)\Phi^r \ud s
- \int_0^t H(s,X^{0,a}_s,{Y}_s,{Z}_s)\Phi_s \ud s \label{eq main error reflexion}
\\
&= & 
\int_0^t \sum_{r=p}^{N_p}\lambda^p_r\set{H_r(s,X^{0,a}_s,Y^r_s,Z^r_s)-H(s,X^{0,a}_s,Y^r_s,Z^r_s)} \Phi^r_s \ud s =:A^p_1 \nonumber
\\ 
&&+ \int_0^t \sum_{r=p}^{N_p}\lambda^p_r \set{H(s,X^{0,a}_s,Y^r_s,Z^r_s)-H(s,X^{0,a}_s,{Y}_s,{Z}_s)} \Phi^r_s \ud s =:A^p_2 \nonumber
\\
&&+ \int_0^t H(s,X^{0,a}_s,{Y}_s,{Z}_s)\set{{}^p\!{\Phi}_s - \Phi_s } \ud s =: A^p_3 \;. \nonumber
\end{eqnarray}
% 
% we have that,
% for all $t \le T$,
% \begin{align}
% \cE^n := \int_0^t  H_n(s,X^{0,a}_s,Y^n_s,Z^n_s)\nabla \varphi_n(Y^n_{s}) \ud s
% &- \int_0^t H(s,X^{0,a}_s,{Y}_s,{Z}_s)\Phi_s \ud s \label{eq main error reflexion}
% \\
% &=
% \\ 
% \int_0^t \set{H_n(s,X^{0,a}_s,Y^n_s,Z^n_s)&-H(s,X^{0,a}_s,Y^n_s,Z^n_s)}\nabla \varphi_n(Y^n_{s}) \ud s =:A^n_1 \nonumber
% \\ 
% + \int_0^t \set{H(s,X^{0,a}_s,Y^n_s,Z^n_s)&-H(s,X^{0,a}_s,{Y}_s,{Z}_s)}\nabla \varphi_n(Y^n_{s}) \ud s =:A^n_2 \nonumber
% \\
% + \int_0^t H(s,X^{0,a}_s,{Y}_s,{Z}_s)&\set{\nabla \varphi_n(Y^n_{s}) - \Phi_s } \ud s =: A^n_3 \;. \nonumber
% \end{align}
We study each term in the right hand side of the above equality separately. For the first one, we compute using Cauchy-Schwartz inequality and the uniform bound on $\NH{2}{\Phi^n}$
\begin{align}
\esp{|A^p_1|} \le C \sum_{r=p}^{N_p}\lambda^p_r\esp{\int_0^t | H_r(s,X^{0,a}_s,Y^r_s,Z^r_s)-H(s,X^{0,a}_s,Y^r_s,Z^r_s) |^2\ud s}^\frac12 \;.
\end{align} 
For all $\kappa>0$, we then get
\begin{align*}
\esp{|A^p_1|} &\le C \sum_{r=p}^{N_p}\lambda^p_r\esp{\int_0^t | H_r(s,X^{0,a}_s,Y^r_s,Z^r_s)-H(s,X^{0,a}_s,Y^r_s,Z^r_s) |^2 \1_{\set{|Y^r_s| + |Z^r_s| \le \kappa}}\ud s}^\frac12 =: B^p_1
\\
&+ C \sum_{r=p}^{N_p}\lambda^p_r\esp{\int_0^t\1_{\set{|Y^r_s| + |Z^r_s| > \kappa}}\ud s}^\frac12 =: B^p_2\,.
\end{align*}
Combining Markov inequality with the uniform square integrability of $Y^n$ and $Z^n$, we easily obtain that 
\begin{align}\label{eq B2 for reflexion}
B^p_2 \le \frac{C}\kappa\;.
\end{align}
For the  term $B^p_1$, we combine the uniform convergence (on compact set) of $H_r$ to $H$, recall Lemma \ref{le smooth f and H}(iii), with the dominated convergence theorem, since $H_r$ and $H$ are bounded, to get that for all $\epsilon>0$ there exists $N_{\kappa,\epsilon} \in \mathbb{N}$ such that 
\begin{align}\label{eq B1 for reflexion}
B^p_1 \le \epsilon \text{ for all } p \ge N_{\kappa,\epsilon}\;.
\end{align}
Combining  \eqref{eq B2 for reflexion} and \eqref{eq B1 for reflexion}, we then get
\begin{align}\label{eq A1 for reflexion}
\lim_p \esp{|A^p_1|} = 0\;.
\end{align}
Next, we compute,  using Cauchy-Schwartz inequality and the uniform bound on $\NH{2}{\Phi^n}$,
\begin{align*}
\esp{|A^p_2|} \le C\sum_{r=p}^{N_p}\lambda^p_r
\esp{\int_0^t |H(s,X^{0,a}_s,Y^r_s,Z^r_s)-H(s,X^{0,a}_s,{Y}_s,{Z}_s) |^2\ud s}^\frac12
\end{align*}
and we deduce
\begin{align}\label{eq A2 for reflexion}
\lim_p \esp{|A^p_2|} = 0\;,
\end{align}
from the continuity of $H$ and the strong convergence of $(Y^r,Z^r)$ to $({Y},{Z})$. 
%The weak convergence of $\nabla \varphi_n (Y^n)$ to $\bar{\Phi}$ allows us to obtain
Finally we use the boundedness of $H$ and the strong convergence of ${}^p\!{\Phi}$ to $\Phi$ to get
\begin{align}\label{eq A3 for reflexion}
\lim_p \esp{|A^p_3|} = 0\;.
\end{align}
Combining  \eqref{eq A1 for reflexion}, \eqref{eq A2 for reflexion} and \eqref{eq A3 for reflexion} 
with \eqref{eq main error reflexion} 
yields $\lim_p \esp{|\cE^p_t|} = 0$.
% which concludes the proof of the characterisation of the reflecting part.
Eventually, we get that, for all $t\le T$,
$$Y_t = g(X_T^{0,a}) + \int_t^T f(s,X_s^{0,a},Y_s,Z_s) \ud s -\int_t^T Z_s \ud W_s - \int_t^T H(s,X_s^{0,a},Y_s,Z_s)\Phi_s \ud s, $$
which concludes the proof of Theorem \ref{theorem existence}. Let us remark that the previous equation allows us to consider a continuous version of the process $Y$.
%converges in $\LP{1}$ to zero as $n \rightarrow \infty$.
%\end{bluetext}
%\textcolor{red}{On prend version continue pour Y ?}
\eproof

\vspace{10pt}
We conclude this section by giving the proof of Corollary \ref{corollaire existence EDSR switching} which is an interesting application of Theorem \ref{theorem existence} to the well studied case of BSDEs for switching problems. Following our approach, the main question reduces now to find an appropriate continuous $H$ to describe the direction of reflection such that $H(\cdot)\Phi = -\Psi$, compare \eqref{BSDE reflected} and \eqref{eqBSDE swithching general}.

\vspace{5pt}
\noindent \textbf{Proof of Corollary \ref{corollaire existence EDSR switching}}\\
It is sufficient to define a continuous function $H$ on $\partial \cD$, recall Remark \ref{re extension H nonverysmooth}. We have
\begin{equation}
\label{def: D construction H}
\cD = \{y \in \R^d : y^l > \max_{j \in \cI} (y^j - c^{lj}), l \in \cI \},
\end{equation}
thus, $\bar{\cD}$ is a non-compact convex polyhedron. We can remark that 
$$\bar{\cD}^0 := \bar{\cD} \cap \{ y^d= 0\}$$
is, by abuse of notation, a compact convex polyhedron of $\R^{d-1}$ and so it is a convex polytope. Indeed, we have
$$\bar{\cD}^0 \subset \{ (y^1,...,y^{d-1}) | y^i \in [-c^{id},c^{di}], \, \forall i \in \{1,...,d-1\} \}\neq \emptyset,$$
since we have $c^{di} + c^{id}>0$ for all $1 \leqslant i \leqslant d-1$ due to the structure condition \eqref{condition de structure c}.
We just have to define $H$ on $\partial \bar{\cD}^0$ and then extend $H$ to $\partial \bar{\cD}$ in this way: for all $(t,x,y,z) \in [0,T] \times \R^k \times \cD \times \R^{d \times k}$, we define
$$H(t,x,y,z) := H(t,x,(y^1-y^d,...,y^{d-1}-y^d,0),z).$$
Since $\bar{\cD}^0$ is a convex polytope, then, by Krein-Milman theorem, it is the convex hull of its extremal points. We will define $H$ on all extremal points and then the value of $H$ on all facets 
$$\mathcal{C}^{lj}=\{y \in \partial \cD^0: y^l = y^j -c^{lj}\},\quad  l,j \in \cI, \quad l\neq j, $$
will be defined by linear interpolations. Let us consider an extremal point $(\bar{y}^1,...,\bar{y}^{d-1})$: we know that there exist $(l_i,j_i)_{i \in \{1,...,d-1\}}\in \{1,...,d\}^{2 \times (d-1)}$ such that 
\begin{itemize}
 \item $(l_i,j_i) \neq (l_k,j_k)$ when $i \neq k$,
 \item for all $i \in \{1,...,d-1\}$, $\bar{y}^{l_i} = \bar{y}^{j_i} -c^{l_i j_i}$ where $\bar{y}^d=0$.
\end{itemize}
Then, we define $H(t,x, (\bar{y}^1,...,\bar{y}^{d-1},0),z)$ as the orthogonal projection onto $\Span (\{e^{l_1},...,e^{l_{d-1}} \})$. To conclude it is sufficient to check that $H(t,x, (y^1,...,y^{d-1},0),z)$ sends the vector $e^l-e^j$ to the vector $e^l$ when $(y^1,...,y^{d-1}) \in \mathcal{C}^{lj}$ and  to show the result only for extremal points. In order to do so, let us consider $(\bar{y}^1,...,\bar{y}^{d-1}) \in \mathcal{C}^{lj}$ an extremal point: by the definition of $H$ we just have to show that $e^j \notin \{e^{l_1},...,e^{l_{d-1}} \}$ where we re-use previous notations. Let us prove it by contradiction: we assume that there exists $i \in \{1,...,d-1\}$ such that
\begin{equation}
\label{eq: construction H 1}
 j=l_i \, \textrm{ and } \, \bar{y}^{l_i} = \bar{y}^{j_i}-c^{l_i j_i}.
\end{equation}
Moreover, we have $(\bar{y}^1,...,\bar{y}^{d-1}) \in \mathcal{C}^{lj}$ so
\begin{equation}
\label{eq: construction H 2}
 \bar{y}^l = \bar{y}^j -c^{lj}.
\end{equation}
By combining \eqref{eq: construction H 1}, \eqref{eq: construction H 2} and the structure condition \eqref{condition de structure c}, we obtain
\begin{eqnarray*}
 \bar{y}^l = \bar{y}^j -c^{lj} = \bar{y}^{j_i} - (c^{lj} + c^{j j_i}) < \bar{y}^{j_i} - c^{lj_i},
\end{eqnarray*}
which is in contradiction with the definition of $\cD$ given by \eqref{def: D construction H}.
\eproof

\subsection{The case of discontinuous $H$}

In this section, we consider the case of a discontinuous direction of reflection on the boundary $\partial D$. We obtain an existence result for an obliquely reflected BSDE but the characterization of the reflecting part is somehow more involved, specially at the discontinuity point of $H$, where many directions of reflection are allowed at the limit. This too weak characterization  leads a to non-uniqueness result as illustrated in the next paragraph. The limiting equation we are studying here is then
\begin{align}
\label{eq BSDE reflected discc H}
 Y_t &= g(X_T^{0,a}) + \int_t^T f(s,X_s^{0,a},Y_s,Z_s) \ud s -\int_t^T Z_s \ud W_s - \int_t^T \Psi_s \ud s, \quad t \in [0,T]\\
 & \Psi_s \in E(s,X_s^{0,a},Y_s,Z_s) \textrm{ and } Y_s \in \bar{\cD} \quad \ud \P \otimes \ud s \text{ a.e.},  \nonumber
\end{align}
with $a \in \R^q$ and, for all $(t,x,y,z) \in [0,T] \times \mathbb{R}^q \times \bar{\cD} \times \mathbb{R}^{d \times k}$,
\begin{align*}
 E(t,x,y,z) := \begin{cases}
               \bigcap_{\varepsilon >0} \pos\left(\left\{H(t,x,\tilde{y},\tilde{z})u | (\tilde{y},\tilde{z} \in B((y,z),\varepsilon), u \in \partial \varphi(y)\right\}\right) &\text{ if } y \in \partial \cD\\
               \{0\} &\text{ if } y \in \cD,
              \end{cases}
\end{align*}
where $\pos(\left\{v_i\right\})$ is the closure of the positive linear span of the family $\left\{v_i\right\}$, and $B(x,\varepsilon)$ is the closed Euclidean ball of center $x$ and radius $\varepsilon$.

\begin{Theorem} 
\label{th:Hdiscontinu}
Assume that assumptions \textbf{(AM)}(i)-(v) hold. Then, there exists a solution in $\SP{2}\times\HP{2}\times\HP{2}$ to \eqref{eq BSDE reflected discc H}.
\end{Theorem}

\begin{Remark}
 When $H$ is continuous, we can easily show that
 \begin{align*}
 E(t,x,y,z) = H(t,x,y,z)\partial \varphi(y)
\end{align*}
which is consistent with the result of Theorem \ref{theorem existence}.
\end{Remark}

%\begin{bluetext}
\proof 
The proof of Theorem \ref{th:Hdiscontinu} strongly follows the proof of Theorem \ref{theorem existence}. The arguments are similar from step 1 to step 5.a. We thus start directly  the proof at step 5.b by studying the reflecting term. 
%Since
%$$  \esp{\int_0^T|H(s,X_s^{0,a},Y^n_s,Z^n_s)\nabla \varphi_n(Y^n_s)|^2 \ud s} \le C,$$
%we have, up to a subsequence, the following weak $L^2([0,T]\times\Omega)$-convergence:
%$$\nabla \varphi_n(Y^n_\cdot ) \rightharpoonup \Psi, \quad \text{when } n \rightarrow + \infty,$$
%and
%$$H(.,X_.^{0,a},Y^n_.,Z^n_.)\nabla \varphi_n(Y^n_\cdot ) \rightharpoonup \Phi, \quad \text{when } n \rightarrow + \infty,$$
%Finally we study the reflecting term.
Since
$$  \esp{\int_0^T|H(s,X_s^{0,a},Y^n_s,Z^n_s)\nabla \varphi_n(Y^n_s)|^2 \ud s} \le C,$$
we have, up to a subsequence, the following weak $L^2([0,T]\times\Omega)$-convergence:
%$$\nabla \varphi_n(Y^n_\cdot ) \rightharpoonup \Psi, \quad \text{when } n \rightarrow + \infty,$$
%and
$$\Psi^n :=H(.,X_.^{0,a},Y^n_.,Z^n_.)\nabla \varphi_n(Y^n_\cdot ) \rightharpoonup \Psi, \quad \text{when } n \rightarrow + \infty.$$
%and we can follow step 2.a in the proof of Proposition \ref{pr smooth and bounded exist orbsde} to obtain
%\begin{align*}
%\Phi \in \partial \varphi(Y) \; \text{ and } \;  \int_0^T  \1_{\set{Y_t \notin \partial \cD}}  |\Phi_t|\ud t = 0 \;,
%\end{align*}
Using once again Mazur's Lemma, we know that there exists a convex combination of $(\Psi^n)_{n \in \mathbb{N}}$ converging strongly in $L^2([0,T]\times\Omega)$, namely
\begin{align*}
 {}^p\!{\Psi} := \sum_{r=p}^{N_p}\lambda^p_r \Psi^r \stackrel{p \rightarrow \infty}{\rightarrow} \Psi,
\end{align*}
where $\lambda^p_r \ge 0$ for all $p \in \mathbb{N}$ and $p \le r \le N_p$, and $\sum_{r=p}^{N_p}\lambda^p_r =1$.
As usual, the following combination
\begin{align*}
 ({}^p\!Y,{}^p\!Z) := \sum_{r=p}^{N_p}\lambda^p_r (Y^r,Z^r)
\end{align*}
still converges to $(Y,Z)$ in $\SP{2} \times \HP{2}$ and, by strong convergence, 
\begin{align*}
 \sum_{r=p}^{N_p}\lambda^p_r f_r(\cdot,X^{0,a},Y^r,Z^r) &\overset{L^1([0,T] \times \Omega, dt \otimes dP)}{\longrightarrow}
  f(\cdot,X^{0,a},Y,Z)\; \text{ and } \; \int_0^t {}^p\!Z_s \ud W_s \overset{\LP{2}}{\longrightarrow} \int_0^t Z_s \ud W_s \;.
\end{align*}
So we can pass to the limit into 
$${}^p\!Y_t = g(X_T^{0,a})+ \int_t^T  \sum_{r=p}^{N_p}\lambda^p_r f_r(s,X^{0,a},Y^r_s,Z^r_s) \ud s -\int_t^T {}^p\!Z_s \ud W_s - \int_t^T {}^p\!\Psi_s \ud s$$
to obtain that
$$Y_t = g(X_T^{0,a}) + \int_t^T f(s,X_s^{0,a},Y_s,Z_s) \ud s -\int_t^T Z_s \ud W_s - \int_t^T \Psi_s ds, \quad \ud t \otimes \ud \P\; \text{a.e.}$$
To conclude we just have to study the direction of reflection. Since we have, for all $n \in \mathbb{N}$,
$$\Psi^n_t :=H(t,X_t^{0,a},Y^n_t,Z^n_t)\nabla \varphi_n(Y^n_t ) \in H(t,X_t^{0,a},Y^n_t,Z^n_t)\partial \varphi(\mathfrak{P}(Y^n_t))$$
and $(Y^n_t,Z^n_t) \stackrel{n \rightarrow \infty}{\rightarrow} (Y,Z)$ $\ud t \otimes \ud \P\; \text{a.e.}$, then, for all $\varepsilon_1>0$ and $\varepsilon_2>0$, there exists $N$ (depending on $\omega$) such that, for all $n \ge N$,
$$\Psi^n_t \in \pos (\{ H(t,X_t^{0,a},\tilde{y},\tilde{z})u | (\tilde{y},\tilde{z}) \in B((Y_t,Z_t),\varepsilon_1), \bar{y} \in B(Y_t,\varepsilon_2) \cap \bar{\cD}, u \in \partial \varphi (\bar{y})\}),$$
$\ud t \otimes \ud \P\; \text{a.e.}$ It implies that, for all $p \ge N$,
$${}^p\!\Psi_t \in \pos (\{ H(t,X_t^{0,a},\tilde{y},\tilde{z})u | (\tilde{y},\tilde{z}) \in B((Y_t,Z_t),\varepsilon_1), \bar{y} \in B(Y_t,\varepsilon_2) \cap \bar{\cD}, u \in \partial \varphi (\bar{y})\})$$
$\ud t \otimes \ud \P\; \text{a.e.}$ Finally we get that 
$$\Psi_t\in \bar{E}(t,X_t^{0,a},Y_t,Z_t) \quad \ud t \otimes \ud \P-\text{a.e.}$$
where
\begin{align*}
 \bar{E}(t,x,y,z) := \bigcap_{\varepsilon_1>0,\varepsilon_2 >0} \pos\left(\left\{H(t,x,\tilde{y},\tilde{z})u | (\tilde{y},\tilde{z} \in B((y,z),\varepsilon_1), \bar{y} \in B(y,\varepsilon_2)\cap\bar{\cD}, u \in \partial \varphi(\bar{y})\right\}\right).
 \end{align*}
 When $y \in \cD$ we can remark that $\partial \varphi(\bar{y})={0}$ when $\bar{y} \in B(y,\varepsilon_2)\cap\bar{\cD}$ with $\varepsilon_2$ small enough: thus we get $E(t,x,y,z)={0}$. {When $y \notin \cD$, Let us show that 
 \begin{equation}
  \label{egalite sous-gradient intersection}
  \partial \varphi(y) = \bigcap_{\varepsilon_2 >0} \left\{ u | u \in \partial \varphi (\bar{y}), \bar{y} \in B(y,\varepsilon_2)\cap\bar{\cD} \right\}.
 \end{equation}
One inclusion is obvious, we will prove the other one. Let us consider $u \in \partial \varphi(y_n)$ for all $n \in \mathbb{N}^*$ with $y_n \in B(y,1/n) \cap \bar{\cD}$ and let us show that $u \in \partial \varphi(y)$. For all $z \in \bar{\cD}$ and  $n \in \mathbb{N}$ we have
\begin{align*}
 u\cdot (z-y) = u\cdot (z-y_n) + u\cdot (y_n-y)
\end{align*}
and so 
\begin{align*} 
 \sup_{z \in \bar{\cD}} \left( u\cdot (z-y) \right) \leqslant \sup_{z \in \bar{\cD}} \left( u\cdot (z-y_n)\right) + u\cdot (y_n-y) \leqslant |u| |y_n-y|
\end{align*}
by definition of $\partial \varphi(y_n)$. Then, by taking $n \rightarrow + \infty$ in the previous inequality we get
\begin{align*} 
 \sup_{z \in \bar{\cD}} \left(u\cdot (z-y)\right) \leqslant 0
\end{align*}
 which proves \eqref{egalite sous-gradient intersection}. This result implies that for any $(\tilde{y},\tilde{z}) \in \R^d \times \R^{d \times k}$,
 \begin{align*}
 &\bigcap_{\varepsilon_2 >0} \left\{H(t,x,\tilde{y},\tilde{z})u | \bar{y} \in B(y,\varepsilon_2)\cap\bar{\cD}, u \in \partial \varphi(\bar{y})\right\}
 = \left\{H(t,x,\tilde{y},\tilde{z})u |  u \in \partial \varphi(y)\right\},
\end{align*}
and so we finally get that $\bar{E} = E$ which concludes the proof.}
\eproof
%\end{bluetext}

\paragraph{A counter-example to uniqueness}

\vspace{4pt}
Inspired by Remark 4.4 in \cite{Lions-Sznitman-84}, we suggest the following counter-example to uniqueness in a non-smooth setting. The domain $\cD$ %(2D view on the $(y_1,y_2)$ plan represented above)
is given by
\begin{align*}
 \cD= \set{y \in \R^3 \,|\, y_1 \ge 0 \text{ and } y_2 + y_1 \ge 0}
\end{align*}
Observe that $\partial \cD = F_1 \cup F_2$, where $F_1$ and $F_2$ are given by
\begin{align*}
 F_1 = \set{y \in \R^3 \,|\, y_1= 0 \text{ and } y_2 \ge 0}
 , \;
 F_2 = \set{y \in \R^3 \,|\, y_1 \ge 0 \text{ and } y_1 + y_2 = 0}
\end{align*}
and we denote by $G= F_1 \cap F_2$, the corner of the domain.
On $F_1$ we assume that the reflection is normal so that $H=I_3$, including points on $G$ where the outward cone of reflection if given by
\begin{align*}
 \cK = \set{ y \in \R^3 \,|\, y_1 \le 0, \, y_2 \le 0 \text{ and } y_2 \ge y_1}\,.
\end{align*}
The direction of reflection is along the $y_1$ axis on $F_2 \setminus G$ and is thus oblique, $H$ is constant but not equal to $I_d$. $H$ is thus discontinuous at the corner. 

\vspace{4pt}
%We consider a BSDE with the following data: $X=W$, $\xi = (0,0,X_T)^\top$, $f(t,x,y,z)=-(1,1,0)^\top$ is constant. Note that it satisfies the assumption of the previous Theorem, and we give now two distinct solutions:
%\begin{enumerate}
 %\item The first solution is given by $Y_t = (0,0,W_t)^\top$ and $K_t = (t,t,0)$.
 %\item The second solution is given by $Y'_t = (t,-t,W_t)^\top$ and $K'_t = (2t,0,0)$.
%\end{enumerate}
%Note that in both cases $Z_t =Z'_t = (0,0,1)$.

\vspace{4pt}
We consider a BSDE with the following data: $X=W$, $\xi = (0,0,X_T)^\top$, $f(t,x,y,z)=-(z^3,z^3,0)^\top$ is constant. Note that it satisfies the assumption \textbf{(AM)}(i)-(v). We give now two distinct solutions:
\begin{enumerate}
 \item The first solution is given by $Y_t = (0,0,W_t)^\top$, $Z_t = (0,0,1)^\top$ and $\Psi_t = (-t,-t,0)^\top$.
 \item The second solution is given by $Y'_t = (T-t,-(T-t),W_t)^\top$, $Z'_t = (0,0,1)^\top$ and $\Psi'_t = (-2t,0,0)^\top$.
\end{enumerate}
%Note that in both cases $Z_t =Z'_t = (0,0,1)$.

%\input{adriensgames}

\bibliographystyle{plain}

\def\cprime{$'$}

\end{document}